\documentclass[a4paper,10pt]{article}%
\usepackage{graphicx}
\usepackage{amsmath}
\usepackage{amsfonts}
\usepackage{amssymb}%
\newtheorem{theorem}{Theorem}

\newtheorem{definition}[theorem]{Definition}
\newtheorem{example}[theorem]{Example}

\newtheorem{lemma}[theorem]{Lemma}

\newtheorem{proposition}[theorem]{Proposition}
\newtheorem{remark}[theorem]{Remark}

\newenvironment{proof}[1][Proof]{\textbf{#1.} }{\ \rule{0.5em}{0.5em}}
\begin{document}

\title{Kernel Theorems in Spaces of Tempered Generalized Functions}
\author{A. Delcroix\\Equipe Analyse Alg\'{e}brique Non Lin\'{e}aire\\\textit{Laboratoire Analyse, Optimisation, Contr\^{o}le}\\Facult\'{e} des sciences - Universit\'{e} des Antilles et de la Guyane\\97159\ Pointe-\`{a}-Pitre Cedex \\Guadeloupe\\E\_mail: Antoine.Delcroix@univ-ag.fr}
\maketitle

\begin{abstract}
In analogy to the classical isomorphism between $\mathcal{L}\left(
\mathcal{S}\left(  \mathbb{R}^{n}\right)  ,\mathcal{S}^{\prime}\left(
\mathbb{R}^{m}\right)  \right)  $ and $\mathcal{S}^{\prime}\left(
\mathbb{R}^{n+m}\right)  $, we show that a large class of moderate linear
mappings acting between the space $\mathcal{G}_{\mathcal{S}}\left(
\mathbb{R}^{n}\right)  $ of Colombeau rapidly decreasing generalized functions
and the space $\mathcal{G}_{\tau}\left(  \mathbb{R}^{n}\right)  $ of temperate
ones admits generalized integral representations, with kernels belonging to
$\mathcal{G}_{\tau}\left(  \mathbb{R}^{n+m}\right)  $. Furthermore, this
result contains the classical one in the sense of the generalized distribution equality.

\end{abstract}

\medskip

\noindent\textbf{Mathematics Subject Classification (2000): 45P05, 46F05,
46F30, 47G10}\smallskip

\noindent\textbf{Keywords:} kernel Theorem, Colombeau temperate generalized
functions, integral operator, temperate distributions.\smallskip

\section{Introduction}

During the three last decades, theories of nonlinear generalized functions
have been developed by many authors (see \cite{AntRad,GKOS,NePiSc,Ober1},...),
mainly based on the ideas of J.-F.\ Colombeau \cite{Col1,Col2}, which we are
going to follow in the sequel. Those theories appear to be a natural
continuation of the distributions' one \cite{HorPDOT1,Schwartz1,Treves1},
specially efficient to pose and solve differential or integral problems with
irregular data. \smallskip

In this paper, we continue the investigations in the field of generalized
integral operators initiated by \cite{Scarpa1} (recently republished in
\cite{Scarpa2,Scarpa3}) and carried on by
\cite{BCD,ADSKTT,Garetto,GaGrOb,Valmo1}. Let us recall that those operators
generalize, in the Colombeau framework, the operators with distributional
kernels in the space of Schwartz distributions \cite{BCD}. More specifically,
in \cite{ADSKTT}, we proved that any moderate net of linear maps $\left(
L_{\varepsilon}:\mathcal{D}\left(  \mathbb{R}^{n}\right)  \rightarrow
\mathrm{C}^{\infty}\left(  \mathbb{R}^{m}\right)  \right)  _{\varepsilon}$,
that is satisfying some growth properties with respect to the parameter
$\varepsilon$, gives rise to a linear map $L:\mathcal{G}_{C}\left(
\mathbb{R}^{n}\right)  \rightarrow\mathcal{G}\left(  \mathbb{R}^{m}\right)  $.
(Where $\mathcal{G}\left(  \mathbb{R}^{d}\right)  $ and $\mathcal{G}%
_{C}\left(  \mathbb{R}^{n}\right)  $\ denote respectively the space of
generalized functions and the space of compactly supported ones.) The main
result is that $L$ can be represented as a generalized integral operator in
the spirit of Schwartz Kernel Theorem. \smallskip

Going \ further in this direction, we study here the generalization of the
classical isomorphism between $\mathcal{S}^{\prime}\left(  \mathbb{R}%
^{n+m}\right)  $ and the space of continuous linear mappings acting between
$\mathcal{S}\left(  \mathbb{R}^{n}\right)  $ and $\mathcal{S}^{\prime}\left(
\mathbb{R}^{m}\right)  $ \cite{Treves1}. Thus, the spaces of generalized
functions considered in this paper are the space $\mathcal{G}_{\mathcal{S}%
}\left(  \mathbb{R}^{n}\right)  $ of rapidly decreasing generalized functions
\cite{ADRapDec,GaGrOb,Radyno} and the space $\mathcal{G}_{\tau}\left(
\mathbb{R}^{n}\right)  $ of tempered generalized functions
\cite{Col2,GKOS,Scarpa1}: $\mathcal{G}_{\mathcal{S}}\left(  \mathbb{R}%
^{n}\right)  $ plays here, roughly speaking, the role of $\mathcal{S}\left(
\mathbb{R}^{n}\right)  $ (resp. $\mathcal{G}_{C}\left(  \mathbb{R}^{n}\right)
$) in the classical case (resp. in \cite{ADSKTT}), whereas $\mathcal{G}_{\tau
}\left(  \mathbb{R}^{n}\right)  $ plays the role of $\mathcal{S}^{\prime
}\left(  \mathbb{R}^{m}\right)  $ (resp. $\mathcal{G}\left(  \mathbb{R}%
^{m}\right)  $) in the classical case (resp. in \cite{ADSKTT}).\smallskip

The main results are the following. First, any kernel $H\in\mathcal{G}_{\tau
}\left(  \mathbb{R}^{m+n}\right)  $ gives rise to a new type of linear
generalized integral operator acting between $\mathcal{G}_{\mathcal{S}}\left(
\mathbb{R}^{n}\right)  $ and $\mathcal{G}_{\tau}\left(  \mathbb{R}^{m}\right)
$ and defined by%
\[
\widetilde{H}:\ \ \mathcal{G}_{\mathcal{S}}\left(  \mathbb{R}^{n}\right)
\rightarrow\mathcal{G}_{\tau}\left(  \mathbb{R}^{m}\right)  ,\;\;f\mapsto
\widetilde{H}(f)=\left[  \left(  x\mapsto\int H_{\varepsilon}%
(x,y)f_{\varepsilon}(y)\,\mathrm{d}y\right)  _{\varepsilon}\right]  _{\tau},
\]
where $\left(  H_{\varepsilon}\right)  _{\varepsilon}$ (resp. $\left(
f_{\varepsilon}\right)  _{\varepsilon}$) is any representative of $H$ (resp.
$f$) and\ $\left[  \cdot\right]  _{\tau}$, the class of an element in
$\mathcal{G}_{\tau}\left(  \mathbb{R}^{m}\right)  $. Moreover, the linear map
$H\mapsto\widetilde{H}$ from $\mathcal{G}_{\tau}\left(  \mathbb{R}%
^{m+n}\right)  $ to the space of linear maps $\mathbf{L}\left(  \mathcal{G}%
_{\mathcal{S}}\left(  \mathbb{R}^{n}\right)  ,\mathcal{G}_{\tau}\left(
\mathbb{R}^{m}\right)  \right)  $ is injective (Proposition \ref{PropOpInt}).
This gives the first part of the expected result. Then, new regular subspaces
of $\mathcal{G}_{\mathcal{S}}\left(  \mathbb{R}^{n}\right)  $ and
$\mathcal{G}_{\tau}\left(  \mathbb{R}^{m}\right)  $ are introduced in the
spirit of \cite{ADRegAppl}. They are used to define the moderate nets of
linear maps, which can be extended to act between $\mathcal{G}_{\mathcal{S}%
}\left(  \mathbb{R}^{n}\right)  $ and $\mathcal{G}_{\tau}\left(
\mathbb{R}^{m}\right)  $ (Proposition \ref{PropExt}).\ Finally, our main
result states that those extensions can be represented as generalized integral
operators (Theorem \ref{KernSTHKern1}). Furthermore, this result is strongly
related to the classical isomorphism theorem recalled above in the following
sense. We can associate to each linear continuous operator $\Lambda
:\mathcal{S}\left(  \mathbb{R}^{n}\right)  \rightarrow\mathcal{S}^{\prime
}\left(  \mathbb{R}^{m}\right)  $, a moderate map $L_{\Lambda}$ and
consequently a kernel $H_{L_{\Lambda}}\in\mathcal{G}_{\tau}\left(
\mathbb{R}^{m+n}\right)  $ such that, for all $f$ in $\mathcal{S}\left(
\mathbb{R}^{n}\right)  $, $\Lambda\left(  f\right)  $ and $\widetilde
{H}_{L_{\Lambda}}\left(  f\right)  $ are equal in the generalized distribution
sense \cite{NePiSc} (Proposition \ref{GSTCoroTHS}).\smallskip

The paper can be divided in two parts.\ The first part, formed by section
\ref{KernSAlgebras} and section \ref{KernsGIO}, introduces the material which
is needed in the sequel (regular spaces of generalized functions, generalized
integral operators). The second part, consisting in the two last sections, is
devoted to the definition of moderate nets, the statement of the main theorems
and their proofs. Concerning them, we insist on the differences with
\cite{ADSKTT}: Replacing $\mathcal{G}_{C}\left(  \mathbb{R}%
^{n}\right)  $ by the bigger space $\mathcal{G}_{\mathcal{S}}\left(
\mathbb{R}^{n}\right)  $ and $\mathcal{G}\left(  \mathbb{R}^{m}\right)  $ by
$\mathcal{G}_{\tau}\left(  \mathbb{R}^{m}\right)  $ obliges to consider global
estimates instead of ones on compact sets.\ This forces to introduce a new concept
of moderate maps, and to refine the estimates and the arguments concerning integration.
However, we also substantially simplify here the proof of the equality in
generalized distribution sense.

\section{Colombeau type algebras\label{KernSAlgebras}}

Throughout this section, $d$ will be a strictly positive integer and $\Omega$
an open subset of $\mathbb{R}^{d}$. As mentioned in the introduction, we only
consider in this paper the spaces $\mathcal{G}_{\mathcal{S}}\left(
\Omega\right)  $ of rapidly decreasing generalized functions and
$\mathcal{G}_{\tau}\left(  \Omega\right)  $ of temperate generalized
functions. For $f\in\mathrm{C}^{\infty}\left(  \Omega\right)  $,
$r\in\mathbb{Z}$ and $l\in\mathbb{N}$, set
\begin{equation}
\mu_{r,l}(f)=\sup_{x\in\Omega,\ \left\vert \alpha\right\vert \leq l}\left(
1+\left\vert x\right\vert \right)  ^{r}\left\vert \partial^{\alpha}f\left(
x\right)  \right\vert \ \ \ \ (\text{with values in }\left[  0,+\infty\right]
). \label{KernSdefSn1}%
\end{equation}

\subsection{Rapidly decreasing generalized functions\label{KernsRDGF}}

Set
\begin{align*}
\mathcal{E}_{\mathcal{S}}\left(  \Omega\right)   &  =\left\{  \left(
f_{\varepsilon}\right)  _{\varepsilon}\in\mathcal{S}\left(  \Omega\right)
^{\left(  0,1\right]  }\,\left\vert \,\forall\left(  q,l\right)  \in
\mathbb{N}^{2},\ \exists N\in\mathbb{N},\;\;\mu_{q,l}\left(  f_{\varepsilon
}\right)  =\mathrm{O}\left(  \varepsilon^{-N}\right)  \;\mathrm{as}%
\;\varepsilon\rightarrow0\right.  \right\} \\
\mathcal{N}_{\mathcal{S}}\left(  \Omega\right)   &  =\left\{  \left(
f_{\varepsilon}\right)  _{\varepsilon}\in\mathcal{S}\left(  \Omega\right)
^{\left(  0,1\right]  }\,\left\vert \,\forall\left(  q,l\right)  \in
\mathbb{N}^{2},\ \forall p\in\mathbb{N},\;\;\mu_{q,l}\left(  f_{\varepsilon
}\right)  =\mathrm{O}\left(  \varepsilon^{p}\right)  \;\mathrm{as}%
\;\varepsilon\rightarrow0\right.  \right\}  .
\end{align*}

One can show that $\mathcal{E}_{\mathcal{S}}\left(  \Omega\right)  $ is a
subalgebra of $\mathcal{S}\left(  \Omega\right)  ^{\left(  0,1\right]  }$ and
$\mathcal{N}_{\mathcal{S}}\left(  \Omega\right)  $ an ideal of $\mathcal{E}%
_{\mathcal{S}}\left(  \Omega\right)  $. The algebra $\mathcal{G}_{\mathcal{S}%
}\left(  \Omega\right)  =\mathcal{E}_{\mathcal{S}}\left(  \Omega\right)
/\mathcal{N}_{\mathcal{S}}\left(  \Omega\right)  $ is called the algebra of
\emph{rapidly decreasing generalized functions} \cite{ADRapDec,GaGrOb,Radyno}.
A straightforward exercise shows that the functor $\mathcal{G}_{\mathcal{S}%
}\left(  \cdot\right)  $ defines a presheaf of differential algebras over
$\mathbb{R}^{d}$ and a presheaf of modules over the factor ring of generalized
constants $\overline{\mathbb{C}}=\mathcal{E}_{M}\left(  \mathbb{C}\right)
/\mathcal{N}\left(  \mathbb{C}\right)  $, with
\begin{align*}
\mathcal{E}_{M}\left(  \mathbb{K}\right)   &  =\left\{  \left(  x_{\varepsilon
}\right)  _{\varepsilon}\in\mathbb{K}^{\left(  0,1\right]  }\,\left\vert
\,\exists N\in\mathbb{N},\ \;\left\vert x_{\varepsilon}\right\vert
=\mathrm{O}\left(  \varepsilon^{-N}\right)  \;\mathrm{as}\;\varepsilon
\rightarrow0\right.  \right\} \\
\mathcal{N}\left(  \mathbb{K}\right)   &  =\left\{  \left(  x_{\varepsilon
}\right)  _{\varepsilon}\in\mathbb{K}^{\left(  0,1\right]  }\,\left\vert
\,\forall p\in\mathbb{N},\ \;\left\vert x_{\varepsilon}\right\vert
=\mathrm{O}\left(  \varepsilon^{p}\right)  \;\mathrm{as}\;\varepsilon
\rightarrow0\right.  \right\}  ,
\end{align*}
for $\mathbb{K}=\mathbb{C}$ or $\mathbb{K}=\mathbb{R},\;\mathbb{R}_{+}%
$.\medskip

Set%
\begin{equation}
\mathcal{N}_{\mathcal{S},\ast}\left(  \Omega\right)  =\left\{  \left(
f_{\varepsilon}\right)  _{\varepsilon}\in\mathrm{C}^{\infty}\left(
\Omega\right)  ^{\left(  0,1\right]  }\,\left\vert \,\forall q\in
\mathbb{N},\;\forall p\in\mathbb{N},\;\;\mu_{q,0}\left(  f_{\varepsilon
}\right)  =\mathrm{O}\left(  \varepsilon^{p}\right)  \;\mathrm{as}%
\;\varepsilon\rightarrow0\right.  \right\}  . \label{GregNSstar}%
\end{equation}
We have the exact analogue of Theorems 1.2.25 and 1.2.27 of \cite{GKOS} (see
\cite{ADRapDec}).

\begin{proposition}
\label{KernGSStar}If the open set $\Omega$ is a box, i.e. the product of $d$
open intervals of $\mathbb{R}$ (bounded or not) then $\mathcal{N}%
_{\mathcal{S}}\left(  \Omega\right)  =\mathcal{N}_{\mathcal{S},\ast}\left(
\Omega\right)  \cap\mathcal{E}_{\mathcal{S}}\left(  \Omega\right)  $.
\end{proposition}

We shall need in the sequel some results concerning embeddings. Consider
$\rho\in\mathcal{S}\left(  \mathbb{R}^{d}\right)  $ which satisfies
\begin{equation}%
{\textstyle\int}
\rho\left(  x\right)  \,\mathrm{d}x=1,\ \ \ \ \
{\textstyle\int}
x^{\alpha}\rho\left(  x\right)  \,\mathrm{d}x=0\;\mathrm{for\;all\;}\alpha
\in\mathbb{N}^{d}\backslash\left\{  0\right\}  \label{CSTGSBinfini}%
\end{equation}
and set
\begin{equation}
\forall\varepsilon\in\left(  0,1\right]  ,\;\;\forall x\in\mathbb{R}%
^{d},\;\;\;\rho_{\varepsilon}\left(  x\right)  =\varepsilon^{-d}\rho\left(
x/\varepsilon^{-1}\right)  . \label{CSTGSBinfini2}%
\end{equation}

\begin{proposition}
\label{KernSPembGS1}\cite{ADRapDec}~\newline(i)~The map\emph{ }%
\[
\sigma_{\mathcal{S}}:\ \ \mathcal{S}\left(  \mathbb{R}^{d}\right)
\rightarrow\mathcal{G}_{\mathcal{S}}\left(  \mathbb{R}^{d}\right)
,\ \ f\mapsto\left[  \left(  f_{\varepsilon}\right)  _{\varepsilon}\right]
_{\mathcal{S}}\ \ \mathrm{with\;}f_{\varepsilon}=f\text{ for all }%
\varepsilon\in\left(  0,1\right]
\]
is an embedding of differential algebras.\newline(ii)~The map
\[
\iota_{\mathcal{S}}:\ \ \mathcal{O}_{C}^{\prime}\left(  \mathbb{R}^{d}\right)
\rightarrow\mathcal{G}_{\mathcal{S}}\left(  \mathbb{R}^{d}\right)
,\ \ u\mapsto\left(  u\ast\rho_{\varepsilon}\right)  _{\varepsilon
}+\mathcal{N}_{\mathcal{S}}\left(  \mathbb{R}^{d}\right)
\]
is an embedding of differential vector spaces. ($\mathcal{O}_{C}^{\prime
}\left(  \mathbb{R}^{d}\right)  $ denotes the space of rapidly decreasing
distributions.)\newline(iii)~Moreover, $\iota_{\mathcal{S}\left\vert
\mathcal{S}\left(  \mathbb{R}^{d}\right)  \right.  }=\sigma_{\mathcal{S}}$.
\end{proposition}

\subsection{Temperate generalized functions}

Set
\begin{gather*}
\mathcal{O}_{M}\left(  \Omega\right)  =\left\{  f\in\mathrm{C}^{\infty}\left(
\Omega\right)  \,\left\vert \,\forall l\in\mathbb{N},\;\exists q\in
\mathbb{N},\;\ \mu_{-q,l}(f)<+\infty\right.  \right\} \\
\mathcal{E}_{\tau}\left(  \Omega\right)  =\left\{  \left(  f_{\varepsilon
}\right)  _{\varepsilon}\in\mathcal{O}_{M}\left(  \Omega\right)  ^{\left(
0,1\right]  }\,\left\vert \,\forall l\in\mathbb{N},\;\exists q\in
\mathbb{N},\;\exists N\in\mathbb{N},\;\ \mu_{-q,l}\left(  f_{\varepsilon
}\right)  =\mathrm{O}\left(  \varepsilon^{-N}\right)  \;\mathrm{as}%
\;\varepsilon\rightarrow0\right.  \right\} \\
\mathcal{N}_{\tau}\left(  \Omega\right)  =\left\{  \left(  f_{\varepsilon
}\right)  _{\varepsilon}\in\mathcal{O}_{M}\left(  \Omega\right)  ^{\left(
0,1\right]  }\,\left\vert \,\forall l\in\mathbb{N},\;\exists q\in
\mathbb{N},\;\forall p\in\mathbb{N},\;\ \mu_{-q,l}\left(  f_{\varepsilon
}\right)  =\mathrm{O}\left(  \varepsilon^{p}\right)  \;\mathrm{as}%
\;\varepsilon\rightarrow0\right.  \right\}  ,
\end{gather*}
where $\left(  \mu_{r,l}\right)  _{\left(  r,l\right)  \in\mathbb{Z}%
\times\mathbb{N}}$ are defined by (\ref{KernSdefSn1}).

The set $\mathcal{O}_{M}\left(  \Omega\right)  $ is the algebra of
multiplicators (or of $\mathrm{C}^{\infty}$ functions with slow growth). One
can show that $\mathcal{E}_{\tau}\left(  \Omega\right)  $ is a subalgebra of
$\mathcal{O}_{M}\left(  \Omega\right)  ^{\left(  0,1\right]  }$ and
$\mathcal{N}_{\tau}\left(  \Omega\right)  $ an ideal of $\mathcal{E}_{\tau
}\left(  \Omega\right)  $. The algebra $\mathcal{G}_{\tau}\left(
\Omega\right)  =\mathcal{E}_{\tau}\left(  \Omega\right)  /\mathcal{N}_{\tau
}\left(  \Omega\right)  $ is called the \emph{algebra of tempered generalized
functions} \cite{Col2,GKOS,Scarpa1}. As for the case of $\mathcal{G}%
_{\mathcal{S}}\left(  \cdot\right)  $, the functor $\mathcal{G}_{\tau}\left(
\cdot\right)  $ defines a presheaf over $\mathbb{R}^{d}$ of differential
algebras and a presheaf of modules over the factor ring $\overline{\mathbb{C}%
}$ introduced in subsection \ref{KernsRDGF}.\medskip

The following result will be useful in the sequel. Set%
\[
\mathcal{N}_{\tau,\ast}\left(  \Omega\right)  =\left\{  \left(  f_{\varepsilon
}\right)  _{\varepsilon}\in\mathcal{O}_{M}\left(  \Omega\right)  ^{\left(
0,1\right]  }\,\left\vert \,\exists q\in\mathbb{N},\;\forall p\in
\mathbb{N},\;\;\mu_{-q,0}\left(  f_{\varepsilon}\right)  =\mathrm{O}\left(
\varepsilon^{p}\right)  \;\mathrm{as}\;\varepsilon\rightarrow0\right.
\right\}  .
\]

\begin{proposition}
\label{KernGtauStar}\cite{GKOS}~If the open set $\Omega$ is a box (see
proposition \ref{KernGSStar})
then $\mathcal{N}_{\tau}\left(  \Omega\right)  =\mathcal{N}_{\tau,\ast}%
\cap\mathcal{E}_{\tau}\left(  \Omega\right)  $.
\end{proposition}

The results about embeddings concerning $\mathcal{G}_{\mathcal{\tau}}\left(
\mathbb{R}^{d}\right)  $ can be summarized (Theorems 1.2.27 and 1.2.28 of
\cite{GKOS}) in the following proposition:

\begin{proposition}
\label{GHTTauEmbed}Consider $\rho\in\mathcal{S}\left(  \mathbb{R}^{d}\right)
$ which satisfies (\ref{CSTGSBinfini}) and define $\left(  \rho_{\varepsilon
}\right)  _{\varepsilon}$ as in (\ref{CSTGSBinfini2}).\newline(i)~The
map\emph{ }%
\[
\sigma_{\mathcal{\tau}}:\ \ \mathcal{O}_{M}\left(  \mathbb{R}^{d}\right)
\rightarrow\mathcal{G}_{\mathcal{\tau}}\left(  \mathbb{R}^{d}\right)
,\ \ f\mapsto\left(  f\right)  _{\varepsilon}+\mathcal{N}_{\mathcal{\tau}%
}\left(  \mathbb{R}^{d}\right)  ,\ \ \mathrm{with\;}f_{\varepsilon}=f\text{
for all }\varepsilon\in\left(  0,1\right]
\]
is an embedding of differential algebras.\newline(ii)~The map
\[
\iota_{\mathcal{\tau}}:\ \ \mathcal{S}^{\prime}\left(  \mathbb{R}^{d}\right)
\rightarrow\mathcal{G}_{\mathcal{\tau}}\left(  \mathbb{R}^{d}\right)
,\ \ u\mapsto\left(  u\ast\rho_{\varepsilon}\right)  _{\varepsilon
}+\mathcal{N}_{\mathcal{\tau}}\left(  \mathbb{R}^{d}\right)
\]
is an embedding of differential vector spaces.\newline(iii)~Moreover,
$\iota_{\tau\left\vert \mathcal{O}_{C}\left(  \mathbb{R}^{d}\right)  \right.
}=\sigma_{\mathcal{\tau}}$, where
\[
\mathcal{O}_{C}\left(  \Omega\right)  =\left\{  f\in\mathrm{C}^{\infty}\left(
\Omega\right)  \,\left\vert \,\exists q\in\mathbb{N},\;\forall l\in
\mathbb{N},\;\ \mu_{-q,l}(f)<+\infty\right.  \right\}
\]

\end{proposition}

\subsection{Regular subalgebras of $\mathcal{G}_{\mathcal{S}}\left(
\Omega\right)  $ and $\mathcal{G}_{\tau}\left(  \Omega\right)  $%
\label{GSTSSPrelST}}

In the sequel, we need to consider some subspaces of $\mathcal{G}%
_{\mathcal{S}}\left(  \cdot\right)  $ and $\mathcal{G}_{\tau}\left(
\cdot\right)  $ with restrictive conditions of growth with respect to
$\varepsilon^{-1}$. These spaces give a good framework for the extension of
linear maps and for the convolution of generalized functions. These are
essential properties for our main result. This notions are new for
$\mathcal{G}_{\tau}\left(  \Omega\right)  $, but the main ideas have been
given (for the purpose of microlocal analysis) for $\mathcal{G}\left(
\cdot\right)  $ and $\mathcal{G}_{\mathcal{S}}\left(  \cdot\right)  $ in
\cite{ADRegAppl,ADRapDec}.

\begin{definition}
\label{GlimDefReg}A non empty subspace $\mathcal{R}$ of $\mathbb{R}%
_{+}^{\mathbb{N}}$ is \emph{regular} if\newline($i)$~$\mathcal{R}$ is
\textquotedblleft overstable\textquotedblright\ by translations and by maximum%
\begin{equation}
\forall N\in\mathcal{R},\ \forall k\in\mathbb{N},\ \exists N^{\prime}%
\in\mathcal{R},\ \forall n\in\mathbb{N},\ \ N\left(  n\right)  +k\leq
N^{\prime}\left(  n\right)  ,\vspace{-0.04in} \label{GlimRegD1a}%
\end{equation}%
\begin{equation}
\forall N_{1}\in\mathcal{R},\ \forall N_{2}\in\mathcal{R},\ \exists
N\in\mathcal{R},\ \forall n\in\mathbb{N},\ \ \max\left(  N_{1}\left(
n\right)  ,N_{2}\left(  n\right)  \right)  \leq N\left(  n\right)  \,;
\label{GlimRegD1b}%
\end{equation}
(ii)~For all $N_{1}\ $and $N_{2}$ in $\mathcal{R}$, there exists
$N\in\mathcal{R}$ such that%
\begin{equation}
\forall\left(  l_{1},l_{2}\right)  \in\mathbb{N}^{2},\ \ N_{1}\left(
l_{1}\right)  +N_{2}\left(  l_{2}\right)  \leq N\left(  l_{1}+l_{2}\right)  .
\label{GlimRegD2}%
\end{equation}

\end{definition}

For example, the set $\mathbb{R}_{+}^{\mathbb{N}}$ of all positive valued
sequences and the set $\mathcal{B}$ of bounded sequences are regular.\medskip

Let $\mathcal{R}$ be a regular subset of $\mathbb{R}_{+}^{\mathbb{N}}$ and
set
\begin{gather*}
\mathcal{E}_{\mathcal{S}}^{\mathcal{R}}\left(  \Omega\right)  =\left\{
\left(  f_{\varepsilon}\right)  \in\mathcal{E}\left(  \Omega\right)
\,\,\left\vert \,\exists N\in\mathcal{R},\ \forall q\in\mathbb{N},\ \forall
l\in\mathbb{N},\;\;\mu_{q,l}\left(  f_{\varepsilon}\right)  =\mathrm{O}\left(
\varepsilon^{-N(l)}\right)  \right.  \right\} \\
\mathcal{E}_{\tau}^{\mathcal{R}}\left(  \Omega\right)  =\left\{  \left(
f_{\varepsilon}\right)  \in\mathcal{E}\left(  \Omega\right)  \,\left\vert
\,\exists N\in\mathcal{R},\ \exists q\in\mathbb{N},\ \forall l\in
\mathbb{N},\;\;\mu_{-q,l}\left(  f_{\varepsilon}\right)  =\mathrm{O}\left(
\varepsilon^{-N(l)}\right)  \right.  \right\}  .
\end{gather*}

\begin{proposition}
\label{GlimPropFReg}~\newline(i)~For all regular subspace $\mathcal{R}$ of
$\mathbb{R}_{+}^{\mathbb{N}}$, $\mathcal{E}_{\mathcal{S}}^{\mathcal{R}}\left(
\cdot\right)  $ (resp. $\mathcal{E}_{\tau}^{\mathcal{R}}\left(  \cdot\right)
$) is a presheaf of differential algebras over the ring $\mathcal{E}%
_{M}\left(  \mathbb{C}\right)  $.\newline(ii)~For all regular subspaces
$\mathcal{R}_{1}$ and $\mathcal{R}_{2}$ of $\mathbb{R}_{+}^{\mathbb{N}}$, with
$\mathcal{R}_{1}\subset\mathcal{R}_{2}$, $\mathcal{E}_{\mathcal{S}%
}^{\mathcal{R}_{1}}\left(  \cdot\right)  $ (resp. $\mathcal{E}_{\tau
}^{\mathcal{R}_{1}}\left(  \cdot\right)  $) is a subpresheaf of $\mathcal{E}%
_{\mathcal{S}}^{\mathcal{R}_{2}}\left(  \cdot\right)  $ (resp. $\mathcal{E}%
_{\tau}^{\mathcal{R}_{2}}\left(  \cdot\right)  $).
\end{proposition}

We refer the reader to \cite{ADRegAppl} for the proof for the case of rapidly
decreasing generalized functions. The proof for temperate ones is similar.

\begin{definition}
\label{GregularAlgebra}For all regular subset $\mathcal{R}$ of $\mathbb{R}%
_{+}^{\mathbb{N}}$, the presheaf of factor algebras $\mathcal{G}_{\mathcal{S}%
}^{\mathcal{R}}\left(  \cdot\right)  $ (resp. $\mathcal{G}_{\tau}%
^{\mathcal{R}}\left(  \cdot\right)  $) is called the presheaf of $\mathcal{R}%
$\emph{-regular algebras of rapidly decreasing }(resp.\emph{ temperate})\emph{
generalized functions}.
\end{definition}

\begin{example}
\label{GlimExample1}Taking $\mathcal{R}=\mathbb{R}_{+}^{\mathbb{N}}$, we
recover the presheaf $\mathcal{G}_{\mathcal{S}}\left(  \cdot\right)  $ (resp.
$\mathcal{G}_{\tau}\left(  \cdot\right)  $) of rapidly decreasing (resp.
temperate) generalized functions. Taking $\mathcal{R}=\mathcal{B}$, we obtain
the presheaf $\mathcal{G}_{\mathcal{S}}^{\infty}\left(  \cdot\right)
$\ (resp. $\mathcal{G}_{\tau}^{\infty}\left(  \cdot\right)  $) which is
analogue for $\mathcal{G}_{\mathcal{S}}\left(  \cdot\right)  $ (resp.
$\mathcal{G}_{\tau}\left(  \cdot\right)  $) to the sheaf of $\mathcal{G}%
^{\infty}$-generalized functions for $\mathcal{G}\left(  \cdot\right)  $. (See
\cite{Ober1} for $\mathcal{G}^{\infty}\left(  \cdot\right)  $ and
\cite{ADRegAppl} for $\mathcal{G}_{\mathcal{S}}^{\infty}\left(  \cdot\right)
$.)
\end{example}

\begin{example}
\textbf{Rapidly decreasing and temperate generalized functions with slow
asymptotic growth}. We consider mainly the two following examples of regular
spaces \cite{ADSKTT}:
\[
\mathcal{L}_{0}=\left\{  q\in\mathbb{R}_{+}^{\mathbb{N}}\mathrm{\;with}\text{
}\lim_{l\rightarrow+\infty}\left(  q(l)/l\right)  =0\right\}
\ \ \ \ \mathcal{L}_{a}=\left\{  q\in\mathbb{R}_{+}^{\mathbb{N}}%
\mathrm{\;with}\text{ }\lim\!\sup_{l\rightarrow+\infty}\left(  q(l)/l\right)
<a\right\}  \ (a>0).
\]
The corresponding presheaves of algebras $\mathcal{G}_{\mathcal{S}%
}^{\mathcal{L}_{a}}\left(  \cdot\right)  $ (resp. $\mathcal{G}_{\tau
}^{\mathcal{L}_{a}}\left(  \cdot\right)  $) are called the presheaves of
rapidly decreasing (resp. temperate) generalized functions with slow
asymptotic growth, with respect to the regularizing parameter $\varepsilon$.
\end{example}

\subsection{Fundamental lemma}

\begin{lemma}
\label{LmnMTh4}Let $a$ be a real in $\left[  0,1\right]  $. Consider $\rho
\in\mathcal{S}\left(  \mathbb{R}^{d}\right)  $ which satisfies
(\ref{CSTGSBinfini}) and define $\left(  \rho_{\varepsilon}\right)
_{\varepsilon}$ as in (\ref{CSTGSBinfini2}).\ For any\textbf{\ }$\left(
g_{\varepsilon}\right)  _{\varepsilon}\in\mathcal{E}_{\mathcal{S}%
}^{\mathcal{L}_{a}}\left(  \mathbb{R}^{d}\right)  $ (resp. \ $\mathcal{E}%
_{\tau}^{\mathcal{L}_{a}}\left(  \mathbb{R}^{d}\right)  $), we have
\begin{equation}
\left(  g_{\varepsilon}\ast\rho_{\varepsilon}-g_{\varepsilon}\right)
_{\varepsilon}\in\mathcal{N}_{\mathcal{S}}\left(  \mathbb{R}^{d}\right)
\ \ (\text{resp. }\mathcal{N}_{\tau}\left(  \mathbb{R}^{d}\right)  \text{).}
\label{GSTConvolEqu}%
\end{equation}

\end{lemma}

\begin{proof}
It suffices to treat the case $a=1$, since $\mathcal{E}_{\mathcal{S}%
}^{\mathcal{L}_{a}}\subset\mathcal{E}_{\mathcal{S}}^{\mathcal{L}_{1}}$ (resp.
$\mathcal{E}_{\tau}^{\mathcal{L}_{a}}\subset\mathcal{E}_{\tau}^{\mathcal{L}%
_{1}}$). We shall do the proof for the case $d=1$, the general case only
differs by more complicated algebraic expressions.\smallskip

\noindent\emph{(i)~Case of }$\mathcal{E}_{\mathcal{S}}^{\mathcal{L}_{1}%
}\left(  \mathbb{R}\right)  $.-\emph{~}Fix $\left(  g_{\varepsilon}\right)
_{\varepsilon}\in\mathcal{E}_{\mathcal{S}}^{\mathcal{L}_{1}}\left(
\mathbb{R}\right)  $ and set, for $\varepsilon\in\left(  0,1\right]  $,
\[
\forall y\in\mathbb{R},\ \ \Delta_{\varepsilon}(y)=\left(  g_{\varepsilon}%
\ast\rho_{\varepsilon}\right)  (y)-g_{\varepsilon}(y)=\int\left(
g_{\varepsilon}(y-x)-g_{\varepsilon}(y)\right)  \rho_{\varepsilon
}(x)\,\mathrm{d}x.
\]
There exists $N\in\mathcal{L}_{1}$ such that%
\[
\forall q\in\mathbb{N},\ \forall l\in\mathbb{N},\ \forall\xi\in\mathbb{R}%
,\ \ \left\vert g_{\varepsilon}^{\left(  l\right)  }\left(  \xi\right)
\right\vert \leq C_{l,q}\left(  1+\left\vert \xi\right\vert \right)
^{-q}\varepsilon^{-N(l)}\text{ \ \ (}C_{l,q}>0),
\]
for $\varepsilon$ smaller than some $\eta_{i,q}$ depending on $i$ and $q$.

Let $p$ and $q$ be two integers. As $\lim\!\sup_{l\rightarrow+\infty}\left(
N\left(  l\right)  /l\right)  <1$, we get $\lim_{l\rightarrow+\infty}\left(
l-N(l)\right)  =+\infty$ and the existence of an integer $k$ such that
$k-N(k)>p$. Taylor's formula gives
\[
\forall\left(  x,y\right)  \in\mathbb{R}^{2},\ \ g_{\varepsilon}%
(y-x)-g_{\varepsilon}(y)=\sum\limits_{l=1}^{k-1}\frac{\left(  -x\right)  ^{l}%
}{l!}g_{\varepsilon}^{\left(  l\right)  }\left(  y\right)  +\frac{\left(
-x\right)  ^{k}}{\left(  k-1\right)  !}\int_{0}^{1}g_{\varepsilon}%
^{(k)}\left(  y-ux\right)  \left(  1-u\right)  ^{k-1}\,\mathrm{d}u,
\]
and
\[
\forall y\in\mathbb{R},\ \ \Delta_{\varepsilon}(y)=\int\frac{\left(
-x\right)  ^{k}}{\left(  k-1\right)  !}(\int_{0}^{1}g_{\varepsilon}%
^{(k)}\left(  y-ux\right)  \left(  1-u\right)  ^{k-1}\,\mathrm{d}%
u)\,\rho_{\varepsilon}(x)\,\mathrm{d}x.
\]
since $\int x^{i}\theta_{\varepsilon}(x)\,\mathrm{d}x=0$ (for all $i\geq1$).
Setting $v=x/\varepsilon$, we get
\[
\forall y\in\mathbb{R},\ \ \Delta_{\varepsilon}(y)=\frac{\varepsilon^{k}%
}{\left(  k-1\right)  !}\int\left(  -v\right)  ^{k}(\int_{0}^{1}%
g_{\varepsilon}^{(k)}\left(  y-\varepsilon uv\right)  \left(  1-u\right)
^{k-1}\,\mathrm{d}u)\,\rho\left(  v\right)  \,\mathrm{d}v.
\]
and%
\begin{equation}
\forall y\in\mathbb{R},\ \ \left\vert \Delta_{\varepsilon}(y)\right\vert
\leq\frac{\varepsilon^{k}}{\left(  k-1\right)  !}\int\left\vert v\right\vert
^{k}(\int_{0}^{1}\left\vert g_{\varepsilon}^{(k)}\left(  y-\varepsilon
uv\right)  \right\vert \,\mathrm{d}u)\,\left\vert \rho\left(  v\right)
\right\vert \,\mathrm{d}v. \label{KernsDemolmFond}%
\end{equation}
For all $y\in\mathbb{R}$, $\varepsilon\in\left(  0,\eta_{k,q}\right]  $,
$u\in\left(  0,1\right]  $, we have
\[
\left\vert g_{\varepsilon}^{(k)}\left(  y-\varepsilon uv\right)  \right\vert
\leq C_{k,q}\left(  1+\left\vert y-\varepsilon uv\right\vert \right)
^{-q}\varepsilon^{-N(k)}.
\]
Since $\rho$ is rapidly decreasing, there exists $C>0$ such that%
\[
\forall v\in\mathbb{R},\ \ \left\vert \rho\left(  v\right)  \right\vert \leq
C\left(  1+\left\vert v\right\vert \right)  ^{-q-k-2}.
\]
Replacing in (\ref{KernsDemolmFond}), we get the existence of a constant
$C_{k,q}^{\prime}>0$ such that
\begin{multline*}
\forall y\in\mathbb{R},\ \forall\varepsilon\in\left(  0,\eta_{k,q}\right]
,\ \ \left\vert \Delta_{\varepsilon}(y)\right\vert \\
\leq C_{k,q}^{\prime}\varepsilon^{k-N(k)}\int(\int_{0}^{1}\left(  1+\left\vert
y-\varepsilon uv\right\vert \right)  ^{-q}\left(  1+\left\vert v\right\vert
\right)  ^{-q}\mathrm{d}u)\,\left\vert v\right\vert ^{k}\left(  1+\left\vert
v\right\vert \right)  ^{-k-2}\,\mathrm{d}v.
\end{multline*}
One can verify (by hand) that $\left(  1+\left\vert y-\varepsilon
uv\right\vert \right)  ^{-q}\left(  1+\left\vert v\right\vert \right)
^{-q}\leq\left(  1+\left\vert y\right\vert \right)  ^{-q}$, for all
$y\in\mathbb{R}$, $\varepsilon\in\left(  0,1\right]  $, $u\in\left(
0,1\right]  $. Thus, there exists a constant $C_{k,q}^{\prime\prime}>0$ such
that%
\[
\forall y\in\mathbb{R},\ \forall\varepsilon\in\left(  0,\eta_{k,q}\right]
,\ \ \left\vert \Delta_{\varepsilon}(y)\right\vert \leq C_{k,q}^{\prime\prime
}\,\varepsilon^{k-N(k)}\left(  1+\left\vert y\right\vert \right)  ^{-q}.
\]
By assumption on $k$, we finally get
\[
\sup_{y\in\mathbb{R}}\left\vert \left(  1+\left\vert y\right\vert \right)
^{q}\Delta_{\varepsilon}(y)\right\vert =\mathrm{O}\left(  \varepsilon
^{p}\right)  \;\mathrm{as}\;\varepsilon\rightarrow0.
\]
Thus, $\Delta_{\varepsilon}(y)$ satisfies the $0$-estimate of the ideal
$\mathcal{N}_{\mathcal{S}}\left(  \mathbb{R}\right)  $. As $\left(
\Delta_{\varepsilon}\right)  _{\varepsilon}\in\mathcal{E}_{\mathcal{S}}\left(
\mathbb{R}\right)  $, we can conclude that $\left(  \Delta_{\varepsilon
}\right)  _{\varepsilon}\in\mathcal{N}_{\mathcal{S}}\left(  \mathbb{R}\right)
$, without estimating the derivatives by using Proposition \ref{KernGSStar}%
.\medskip

\noindent\emph{(ii)~Case of }$\mathcal{E}_{\tau}^{\mathcal{L}_{1}}\left(
\mathbb{R}\right)  $.-\emph{~}The proof is an improvement of the proof of
Theorems 1.2.28 of \cite{GKOS}, based on the ideas developed for the case of
$\mathcal{E}_{\emph{S}}^{\mathcal{L}_{1}}\left(  \mathbb{R}\right)  $ above.
\end{proof}

\begin{remark}
\label{GSTRemKernIdent}Consider a net of mollifiers $\left(  \rho
_{\varepsilon}\right)  _{\varepsilon}$ as in Lemma \ref{LmnMTh4}.\ Relation
(\ref{GSTConvolEqu}) shows that $\left[  \left(  \rho_{\varepsilon}\right)
_{\varepsilon}\right]  _{\mathcal{S}}$ \ (resp. $\left[  \left(
\rho_{\varepsilon}\right)  _{\varepsilon}\right]  _{\tau}$) plays the role of
identity for convolution in $\mathcal{G}_{\mathcal{S}}^{\mathcal{L}_{a}%
}\left(  \mathbb{R}^{d}\right)  $ (resp. $\mathcal{G}_{\tau}^{\mathcal{L}_{a}%
}\left(  \mathbb{R}^{d}\right)  $) whereas this is false for $\mathcal{G}%
_{\mathcal{S}}\left(  \mathbb{R}^{d}\right)  $ (resp. $\mathcal{G}_{\tau
}\left(  \mathbb{R}^{d}\right)  $). This is an essential feature of these new spaces.
\end{remark}

\section{Generalized integral operators\label{KernsGIO}}

As mentioned in the introduction, we consider here generalized integral
operators acting on $\mathcal{G}_{\mathcal{S}}(\mathbb{R}^{n})$ with values in
$\mathcal{G}_{\tau}(\mathbb{R}^{m})$. (We refer the reader to \cite{BCD} and
\cite{GaGrOb} for the more usual case of generalized integral operators acting
on $\mathcal{G}(\mathbb{R}^{n})$.) From now on $m$ and $n$ are two strictly
positive integers.

\begin{lemma}
\label{KernSLmmOps1}Consider $H\in\mathcal{G}_{\tau}(\mathbb{R}^{m+n})$,
$f\in\mathcal{G}_{\mathcal{S}}(\mathbb{R}^{n})$ and $\left(  H_{\varepsilon
}\right)  _{\varepsilon}$ (resp. $\left(  f_{\varepsilon}\right)
_{\varepsilon}$) any representative of $H$ (resp. $f$). The net of
\textrm{C}$^{\infty}$ maps
\begin{equation}
\left(  \widetilde{H}_{\varepsilon}\left(  f_{\varepsilon}\right)  \right)
_{\varepsilon}:=\left(  x\mapsto\int H_{\varepsilon}(x,y)f_{\varepsilon
}(y)\,\mathrm{d}y\right)  _{\varepsilon} \label{KernSLmmOps2}%
\end{equation}
belongs to $\mathcal{E}_{\tau}(\mathbb{R}^{n})$ and the class $\left[  \left(
\widetilde{H}_{\varepsilon}\left(  f_{\varepsilon}\right)  \right)
_{\varepsilon}\right]  _{\tau}$ depends only on $H$ and $f$ but not on the
representatives $\left(  H_{\varepsilon}\right)  _{\varepsilon}$ and $\left(
f_{\varepsilon}\right)  _{\varepsilon}$.
\end{lemma}

\begin{proof}
Firstly, for all $x\in\mathbb{R}^{m}$ and $\varepsilon\in\left(  0,1\right]
$, $H_{\varepsilon}(x,\cdot)f_{\varepsilon}(\cdot)$ belongs to $L^{1}\left(
\mathbb{R}^{n}\right)  $ as well as its derivatives, since $f_{\varepsilon}%
\in\mathcal{S}(\mathbb{R}^{n})$ and $H_{\varepsilon}(x,\cdot)\in
\mathcal{O}_{M}(\mathbb{R}^{n})$. Thus, $\widetilde{H}_{\varepsilon}\left(
f_{\varepsilon}\right)  $ is well defined and it is easily seen that
$\widetilde{H}_{\varepsilon}\left(  f_{\varepsilon}\right)  $ belongs to
$\mathrm{C}^{\infty}(\mathbb{R}^{m})$. Secondly, for any $\alpha\in
\mathbb{N}^{m}$, there exist $q_{1}\in\mathbb{N}$ and $C_{1}>0$ such that (for
$\varepsilon$ small enough)%
\begin{equation}
\forall\left(  x,y\right)  \in\mathbb{R}^{m+n},\text{\ \ \ }\left\vert
\partial_{x}^{\alpha}H_{\varepsilon}(x,y)\right\vert \leq C_{1}\,\left(
1+\left\vert x\right\vert \right)  ^{q_{1}}\left(  1+\left\vert y\right\vert
\right)  ^{q_{1}}\varepsilon^{-q_{1}}. \label{KernSlmnOpsP1}%
\end{equation}
As $\left(  f_{\varepsilon}\right)  _{\varepsilon}\in\mathcal{E}_{\mathcal{S}%
}(\mathbb{R}^{n})$, there exist $q_{2}\in\mathbb{N}$ and $C_{2}>0$ such that
(for $\varepsilon$ small enough)%
\begin{equation}
\forall y\in\mathbb{R}^{n},\text{\ \ \ }\left\vert f_{\varepsilon
}(y)\right\vert \leq C_{2}\,\left(  1+\left\vert y\right\vert \right)
^{-q_{1}-d-1}\varepsilon^{-q_{2}}. \label{KernSlmnOpsP2}%
\end{equation}
Inserting (\ref{KernSlmnOpsP1}) and (\ref{KernSlmnOpsP2}) in the definition of
$\widetilde{H}_{\varepsilon}\left(  f_{\varepsilon}\right)  $, we get a
constant $C_{3}>0$ such that
\[
\forall x\in\mathbb{R}^{m},\text{\ \ \ }\left\vert \partial^{\alpha}%
\widetilde{H}_{\varepsilon}\left(  f_{\varepsilon}\right)  (x)\right\vert \leq
C_{3}\,\left(  1+\left\vert x\right\vert \right)  ^{q_{1}}\varepsilon
^{-q_{1}-q_{2}}.
\]
Thus $\left(  \widetilde{H}_{\varepsilon}\right)  _{\varepsilon}$ is in
$\mathcal{E}_{\tau}(\mathbb{R}^{m})$. Finally, suppose that $\left(
H_{\varepsilon}\right)  _{\varepsilon}$ (resp. $\left(  f_{\varepsilon
}\right)  _{\varepsilon}$) is in $\mathcal{N}_{\tau}(\mathbb{R}^{m+n})$ (resp.
$\mathcal{N}_{\mathcal{S}}(\mathbb{R}^{n})$).\ From estimates similar to
(\ref{KernSlmnOpsP1}) and (\ref{KernSlmnOpsP2}), we get that $\left(
\widetilde{H}_{\varepsilon}\left(  f_{\varepsilon}\right)  \right)
_{\varepsilon}$ is in $\mathcal{N}_{\tau}(\mathbb{R}^{m})$.\ Hence, the
independence on the representatives of $\left[  \left(  \widetilde
{H}_{\varepsilon}\left(  f_{\varepsilon}\right)  \right)  _{\varepsilon
}\right]  _{\tau}$ is proved.
\end{proof}

\medskip Lemma\ \ref{KernSLmmOps1} justifies the following:

\begin{definition}
\label{KernDefOps1}Let $H$ be in $\mathcal{G}_{\tau}(\mathbb{R}^{m+n})$. The
integral operator of kernel $H$ is the map $\widetilde{H}$ defined by
\[
\widetilde{H}:\ \ \mathcal{G}_{\mathcal{S}}\left(  \mathbb{R}^{n}\right)
\rightarrow\mathcal{G}_{\tau}\left(  \mathbb{R}^{m}\right)  ,\ \ \ f\mapsto
\widetilde{H}\left(  f\right)  =\left[  \left(  x\mapsto\int H_{\varepsilon
}(x,y)f_{\varepsilon}(y)\,\mathrm{d}y\right)  _{\varepsilon}\right]  _{\tau},
\]
with the notations of Lemma \ref{KernSLmmOps1}.
\end{definition}

\begin{proposition}
\label{PropOpInt}With the notations of Definition \ref{KernDefOps1}, the
operator $\widetilde{H}$ defines a linear mapping from $\mathcal{G}_{\tau
}\left(  \mathbb{R}^{n}\right)  $ to $\mathcal{G}_{\mathcal{S}}\left(
\mathbb{R}^{m}\right)  $. Moreover, the map
\[
\mathcal{G}_{\tau}(\mathbb{R}^{m+n})\rightarrow\mathbf{L}\left(
\mathcal{G}_{\mathcal{S}}\left(  \mathbb{R}^{n}\right)  ,\mathcal{G}_{\tau
}\left(  \mathbb{R}^{m}\right)  \right)  ,\;\;H\mapsto\widetilde{H}%
\]
is injective.
\end{proposition}

\begin{proof}
Only the last assertion needs a proof. Consider $H\in\mathcal{G}_{\tau
}(\mathbb{R}^{m+n})$ and let $\left(  H_{\varepsilon}\right)  _{\varepsilon}$
be one of its representative. As $\left(  H_{\varepsilon}\right)
_{\varepsilon}$ is in $\mathcal{E}_{\tau}(\mathbb{R}^{m+n})$, there exist
$q\in\mathbb{N}$ and $C_{0}>0$ such that%
\[
\forall\left(  x,y\right)  \in\mathbb{R}^{m+n},\ \ \left\vert H_{\varepsilon
}\left(  x,y\right)  \right\vert \leq C_{0}\left(  1+\left\vert \left(
x,y\right)  \right\vert ^{2}\right)  ^{q}\varepsilon^{-q}\leq C_{0}\left(
\left(  1+\left\vert x\right\vert ^{2}\right)  \left(  1+\left\vert
y\right\vert ^{2}\right)  \right)  ^{q}\varepsilon^{-q},
\]
for $\varepsilon$ small enough. Now, suppose that $\widetilde{H}$ is null.
Then, for all $g\in\mathcal{G}_{\mathcal{S}}\left(  \mathbb{R}^{n}\right)  $
with representative $\left(  g_{\varepsilon}\right)  _{\varepsilon}$, there
exists $r>0$ such that%
\begin{equation}
\forall p\in\mathbb{N},\ \exists C>0,\ \forall x\in\mathbb{R}^{m}%
,\ \ \left\vert \widetilde{H}_{\varepsilon}\left(  g_{\varepsilon}\right)
\left(  x\right)  \right\vert =\left\vert \int H_{\varepsilon}%
(x,y)g_{\varepsilon}(y)\,\mathrm{d}y\right\vert \leq C\left(  1+\left\vert
x\right\vert ^{2}\right)  ^{r}\varepsilon^{p}, \label{KernsPrInj1}%
\end{equation}
for $\varepsilon$ small enough.\smallskip

Set, for all $\left(  x,y\right)  \in\mathbb{R}^{m+n}$, $B_{\varepsilon
}(x,y)=\varepsilon^{q}H_{\varepsilon}(x,y)\left(  \left(  1+\left\vert
x\right\vert ^{2}\right)  \left(  1+\left\vert y\right\vert ^{2}\right)
\right)  ^{-q-1}$. The net $\left(  B_{\varepsilon}\right)  _{\varepsilon}$
belongs to $\mathcal{E}_{\tau}(\mathbb{R}^{m+n})$. For $f\in\mathcal{G}%
_{\mathcal{S}}\left(  \mathbb{R}^{n}\right)  $ with representative $\left(
f_{\varepsilon}\right)  _{\varepsilon}$, we have
\[
\int B_{\varepsilon}(x,y)f_{\varepsilon}(y)\,\mathrm{d}y=\left(  1+\left\vert
x\right\vert ^{2}\right)  ^{-q-1}\int H_{\varepsilon}(x,y)\,\varepsilon
^{q}\left(  1+\left\vert y\right\vert ^{2}\right)  ^{-q-1}f_{\varepsilon
}(y)\,\mathrm{d}y.
\]
As the net of functions $\left(  y\mapsto\varepsilon^{q}\left(  1+\left\vert
y\right\vert ^{2}\right)  ^{-q-1}f_{\varepsilon}(y)\right)  _{\varepsilon}$
belongs to $\mathcal{E}_{\mathcal{S}}\left(  \mathbb{R}^{n}\right)  $, we get
from (\ref{KernsPrInj1}), the existence of $r_{1}\in\mathbb{Z}$ ($r_{1}%
=r-q-1$) such that
\begin{equation}
\forall p\in\mathbb{N},\ \exists C>0,\ \forall x\in\mathbb{R}^{m}%
,\ \ \left\vert \widetilde{B}_{\varepsilon}\left(  f_{\varepsilon}\right)
\left(  x\right)  \right\vert =\left\vert \int B_{\varepsilon}%
(x,y)f_{\varepsilon}(y)\,\mathrm{d}y\right\vert \leq C\left(  1+\left\vert
x\right\vert ^{2}\right)  ^{r_{1}}\varepsilon^{p}, \label{KernsPrInj1bis}%
\end{equation}
for $\varepsilon$ small enough. Thus, the operator $\widetilde{B}$ defined by
$B=\left[  \left(  B_{\varepsilon}\right)  _{\varepsilon}\right]  _{\tau}$ is null.\ 

Remark that it suffices to show that
\begin{equation}
\forall p\in\mathbb{N},\ \exists\eta^{\prime}>0,\ \exists C^{\prime
}>0,\ \forall\varepsilon\in\left(  0,\eta^{\prime}\right]  ,\ \forall\left(
x,y\right)  \in\mathbb{R}^{m+n},\ \ \ \left\vert B_{\varepsilon}\left(
x,y\right)  \right\vert \leq C^{\prime}\varepsilon^{p} \label{KernsPrInj2}%
\end{equation}
to conclude that $H$ is null. Indeed, if (\ref{KernsPrInj2}) holds, for any
$p\in\mathbb{N}$, we get a constant $C^{\prime}>0$ such that
\[
\forall\left(  x,y\right)  \in\mathbb{R}^{m+n},\ \left\vert H_{\varepsilon
}\left(  x,y\right)  \right\vert \leq\varepsilon^{p-q}C^{\prime}\left(
\left(  1+\left\vert x\right\vert ^{2}\right)  \left(  1+\left\vert
y\right\vert ^{2}\right)  \right)  ^{q+1}\leq\varepsilon^{p-q}C^{\prime
}\left(  1+\left\vert \left(  x,y\right)  \right\vert ^{2}\right)  ^{2q+2},
\]
for $\varepsilon$ small enough. As $q$ is fixed, this proves that $\left(
H_{\varepsilon}\right)  _{\varepsilon}$ satisfies the $0$-estimate of
$\mathcal{N}_{\tau}(\mathbb{R}^{m+n})$. Furthermore, $\left(  H_{\varepsilon
}\right)  _{\varepsilon}$ is in $\mathcal{E}_{\tau}(\mathbb{R}^{m+n})$: Thus,
Proposition \ref{KernGtauStar} shows that $\left(  H_{\varepsilon}\right)
_{\varepsilon}$ is in $\mathcal{N}_{\tau}(\mathbb{R}^{m+n})$.\smallskip
\smallskip

We are now going to show (\ref{KernsPrInj2}). We suppose, on the contrary,
that
\begin{equation}
\exists p_{0}\in\mathbb{N},\ \forall\eta^{\prime}>0,\ \forall C^{\prime
}>0,\ \exists\varepsilon\in\left(  0,\eta^{\prime}\right]  ,\ \exists\left(
x,y\right)  \in\mathbb{R}^{m+n},\ \ \ \left\vert B_{\varepsilon}\left(
x,y\right)  \right\vert >C^{\prime}\varepsilon^{p_{0}}. \label{KernsPrinj30}%
\end{equation}
Thus, there exists a sequence $\left(  \varepsilon_{l}\right)  _{l}$
converging to $0$ such that, for all $l\in\mathbb{N}$, there exists $\left(
x_{l},y_{l}\right)  \in\mathbb{R}^{m+n}$ such that $\left\vert B_{\varepsilon
_{l}}\left(  x_{l},y_{l}\right)  \right\vert \geq\varepsilon_{l}^{p_{0}}$.

As
\begin{equation}
\forall\varepsilon\in\left(  0,1\right]  ,\ \ \forall\left(  x,y\right)
\in\mathbb{R}^{m+n},\ \ \ \left\vert B_{\varepsilon}\left(  x,y\right)
\right\vert \leq C_{0}\left(  1+\left\vert x\right\vert ^{2}\right)
^{-1}\left(  1+\left\vert y\right\vert ^{2}\right)  ^{-1}, \label{KernsPrInj3}%
\end{equation}
the function $B_{\varepsilon}$\ is bounded by $C_{0}$, with $\lim_{\left\vert
\left(  x,y\right)  \right\vert \rightarrow+\infty}\left\vert B_{\varepsilon
}(x,y)\right\vert =0$.\ Thus, there exists $\left(  x_{\varepsilon
},y_{\varepsilon}\right)  \in\mathbb{R}^{m+n}$, such that $\sup_{\left(
x,y\right)  \in\mathbb{R}^{m+n}}\left\vert B_{\varepsilon}(x,y)\right\vert
=\left\vert B_{\varepsilon}\left(  x_{\varepsilon},y_{\varepsilon}\right)
\right\vert =M_{\varepsilon}$. Moreover, due to (\ref{KernsPrinj30}), we have%
\begin{equation}
\forall l\in\mathbb{N},\ \ \ M_{\varepsilon_{l}}>\varepsilon_{l}^{p_{0}}.
\label{KernsPrInj3a}%
\end{equation}

Note that we necessarily have
\begin{equation}
\forall l\in\mathbb{N},\ \ \left\vert x_{\varepsilon_{l}}\right\vert
\leq\varepsilon^{-p_{0}}\,\ \mathrm{and}\ \ \left\vert y_{\varepsilon_{l}%
}\right\vert \leq\varepsilon^{-p_{0}}. \label{KernsPrInj3b}%
\end{equation}
(Indeed, supposing that: $\exists l\in\mathbb{N},\ $ $\left\vert
x_{\varepsilon_{l}}\right\vert >\varepsilon_{l}^{-p_{0}}\,$or $\left\vert
y_{\varepsilon}\right\vert >\varepsilon_{l}^{-p_{0}}$ contradicts
(\ref{KernsPrInj3a}), using (\ref{KernsPrInj3}).)

For all $l\in\mathbb{N}$, we can find a neighborhood $V_{l}$ (resp. $W_{l}$)
of $x_{\varepsilon_{l}}$ (resp. $y_{\varepsilon_{l}}$) such that
\[
\forall\left(  x,y\right)  \in V_{l}\times W_{l},\ \ \ \left\vert
B_{\varepsilon}(x,y)\right\vert \geq M_{\varepsilon_{l}}/2.
\]
Moreover, for all $l\in\mathbb{N}$, $V_{l}$ (resp. $W_{l}$) can be chosen such
that its diameter $\delta\left(  W_{l}\right)  $ is greater than some
$\varepsilon_{l}^{p_{1}}$, for a fixed $p_{1}\geq0$. Indeed, as the net
$\left(  B_{\varepsilon}\right)  _{\varepsilon}$ belongs to $\mathcal{E}%
_{\tau}(\mathbb{R}^{m+n})$ and as $\left(  x_{\varepsilon_{l}},y_{\varepsilon
_{l}}\right)  _{\varepsilon}$ satisfies (\ref{KernsPrInj3b}), the differential
of $B_{\varepsilon_{l}}$ grows at most polynomially in $\varepsilon_{l}^{-1}$
for $l\rightarrow+\infty$ in some convex neighborhood of $\left(
x_{\varepsilon_{l}},y_{\varepsilon_{l}}\right)  _{\varepsilon}$ of diameter,
let us say, $1$. Thus, supposing that: $\forall s\in\mathbb{N},\ \exists
l\in\mathbb{N},\ $\ $\delta\left(  W_{l}\right)  \leq\varepsilon_{l}^{s}$
leads to a contradiction, since it violates the growth property of the
sequence $\left(  \nabla B_{\varepsilon_{l}}\right)  _{l}$.

We define a net $\left(  \theta_{\varepsilon}\right)  _{\varepsilon}$ as
follows. For all $l\in\mathbb{N}$, $\theta_{\varepsilon_{l}}\in\mathcal{D}%
\left(  \mathbb{R}^{n}\right)  $ is such that $0\leq\theta_{\varepsilon_{l}%
}\leq1$, $\theta_{\varepsilon}\equiv1$ on $W_{l}$. For $\varepsilon
\notin\left\{  \varepsilon_{l},l\in\mathbb{N}\right\}  $ we simply choose
$\theta_{\varepsilon}\equiv0$.\ Moreover, we can suppose that the net $\left(
\theta_{\varepsilon}\right)  _{\varepsilon}$ is in $\mathcal{E}_{\mathcal{S}%
}\left(  \mathbb{R}^{n}\right)  $. (This is done by starting from some
$\theta\in\mathcal{D}\left(  \mathbb{R}^{n}\right)  $ with $0\leq\theta\leq1$,
$\theta\equiv1$ on $B(0,1)$ and then by using linear transformations.) Set,
for all $\varepsilon\in\left(  0,1\right]  $ and $y\in\mathbb{R}^{n}$,
$f_{\varepsilon}(y)=\overline{B_{\varepsilon}(x_{\varepsilon},y)}%
\theta_{\varepsilon}\left(  y\right)  $. As the net $\left(  B_{\varepsilon
}\right)  _{\varepsilon}$ belongs to $\mathcal{E}_{\tau}(\mathbb{R}^{n})$ and
$\left(  x_{\varepsilon_{l}}\right)  _{l}$ satisfies (\ref{KernsPrInj3b}), the
net $\left(  f_{\varepsilon}\right)  _{\varepsilon}$ is in $\mathcal{E}%
_{\mathcal{S}}(\mathbb{R}^{n})$.\ Therefore, using (\ref{KernsPrInj1bis}), we
get the existence of $r_{1}$ -\thinspace which can always be supposed
positive\thinspace- such that
\begin{equation}
\forall p\in\mathbb{N},\ \exists\eta>0,\ \exists C>0,\ \forall x\in
\mathbb{R}^{m},\ \forall\varepsilon\in\left(  0,\eta\right]
,\ \ \ \ \left\vert \widetilde{B}_{\varepsilon}\left(  f_{\varepsilon}\right)
\left(  x\right)  \right\vert \leq C\left(  1+\left\vert x\right\vert
^{2}\right)  ^{r_{1}}\varepsilon^{p}, \label{KernsPrInj5}%
\end{equation}
Returning to the definition of $\left(  W_{l}\right)  _{l}$and $\left(
\theta_{\varepsilon_{l}}\right)  _{l}$, we have, for all $l\in\mathbb{N}$,%
\begin{align}
\left.  \left\vert \widetilde{B}_{\varepsilon_{l}}\left(  f_{\varepsilon_{l}%
}\right)  \left(  x_{\varepsilon_{l}}\right)  \right\vert =\int\left\vert
B_{\varepsilon_{l}}(x_{\varepsilon_{l}},y)\right\vert ^{2}\theta
_{\varepsilon_{l}}\left(  y\right)  \,\mathrm{d}y\right.   &  \geq\int_{W_{l}%
}\left\vert B_{\varepsilon_{l}}(x_{\varepsilon_{l}},y)\right\vert
^{2}\,\mathrm{d}y\nonumber\\
&  \geq\delta\left(  W_{l}\right)  M_{\varepsilon_{l}}^{2}/4\geq
\varepsilon_{l}^{p_{1}}M_{\varepsilon_{l}}^{2}/4. \label{KernsPrInj7}%
\end{align}
Inserting (\ref{KernsPrInj7}) in (\ref{KernsPrInj5}), we get
\[
\forall p\in\mathbb{N},\ \ \exists C>0,\ \ M_{\varepsilon_{l}}^{2}%
\leq4C\left(  1+\left\vert x_{\varepsilon_{l}}\right\vert ^{2}\right)
^{r_{1}}\varepsilon_{l}^{p-p_{1}},
\]
for $l$ large enough. From (\ref{KernsPrInj3b}), we have $\left(  1+\left\vert
x_{\varepsilon_{l}}\right\vert ^{2}\right)  ^{r_{1}}\leq\left(  1+\varepsilon
_{l}^{-2p_{0}}\right)  ^{r_{1}}\leq2^{r_{1}}\varepsilon_{l}^{-2p_{0}r_{1}}$,
for all $l\in\mathbb{N}$.\ Finally, by setting $p_{2}=2p_{0}r_{1}+p_{1}$, we
get
\[
\forall p\in\mathbb{N},\ \ \exists C^{\prime}>0,\ \ M_{\varepsilon_{l}}\leq
C^{\prime}\varepsilon_{l}^{\left(  p-p_{2}\right)  /2},
\]
for $l$ large enough. Thus, we get a contradiction with (\ref{KernsPrInj3a}),
ending the proof.
\end{proof}

\section{Kernel Theorems\label{KernsKernt}}

\subsection{Extension of linear maps\label{SSBExtension}}

Nets of maps $\left(  L_{\varepsilon}\right)  _{\varepsilon}$ between two
topological algebras, having some good growth properties with respect to the
parameter $\varepsilon$, can be extended to act between the corresponding
Colombeau algebras, as it is shown in \cite{ADSKTT,GKOS,Scarpa1} for example.
We are going to introduce here new notions adapted to our framework. In the
sequel $\mathcal{L}\left(  \cdot,\cdot\right)  $ (resp. $\mathbf{L}\left(
\cdot,\cdot\right)  $) denote a space of continuous linear maps acting between
classical spaces (resp. of $\overline{\mathbb{C}}$ linear maps acting between
generalized spaces).

\begin{definition}
\label{DefExt1}Let $j$ be an integer and $\left(  L_{\varepsilon}\right)
_{\varepsilon}\in\mathcal{L}\left(  \mathcal{S}\left(  \mathbb{R}^{n}\right)
,\mathcal{O}_{M}\left(  \mathbb{R}^{m}\right)  \right)  ^{\left(  0,1\right]
}$ be a net of linear maps.\newline$\left(  i\right)  $~We say that $\left(
L_{\varepsilon}\right)  _{\varepsilon}$ is\emph{\ moderate}
(resp.\ \emph{negligible}) if%
\begin{equation}%
\begin{array}
[c]{c}%
\forall l\in\mathbb{N},\;\exists\left(  C_{\varepsilon}\right)  _{\varepsilon
}\in\mathcal{E}_{M}\left(  \mathbb{R}_{+}\right)  \ \ \text{(resp.\ }%
\mathcal{N}\left(  \mathbb{R}_{+}\right)  \text{),\ }\exists\left(
p,q,l^{\prime}\right)  \in\mathbb{N}^{3},~~~~~~~~~~~~~~~~~~\\
~~~~~~~~~~~~~~~~~~~~~\forall f\in\mathcal{S}\left(  \mathbb{R}^{n}\right)
,\;\;\mu_{-p,l}\left(  L_{\varepsilon}\left(  f\right)  \right)  \leq
C_{\varepsilon}\,\mu_{q,l^{\prime}}\left(  f\right)  \text{,}\,\text{ for
}\varepsilon\text{ small enough}.
\end{array}
\label{Defmod1}%
\end{equation}
$\left(  ii\right)  $~Let $\left(  b,c\right)  $ be in $\left[  0,+\infty
\right]  \times\mathbb{R}_{+}$. We say that $\left(  L_{\varepsilon}\right)
_{\varepsilon}$ is $\mathcal{L}_{b,c}$-\emph{strongly moderate }if%
\begin{equation}%
\begin{array}
[c]{c}%
\exists\lambda\in\mathcal{L}_{b}\,,\;\exists r\in\mathcal{L}_{c}\,,\;\forall
l\in\mathbb{N},\;\exists C\in\mathbb{R}_{+},\ \ \exists\left(  p,q\right)
\in\mathbb{N}^{2},~~~~~~~~~~~~~~~~~~~~~~~~~~\\
~~~~~~~~~~~~~~\forall f\in\mathcal{S}\left(  \mathbb{R}^{n}\right)
,\;\;\mu_{-p,l}\left(  L_{\varepsilon}\left(  f\right)  \right)  \leq
C\,\varepsilon^{-r(l)}\,\mu_{q,\lambda\left(  l\right)  }\left(  f\right)
\,\text{,}\,\text{ for }\varepsilon\text{ small enough}.
\end{array}
\label{Defmod2}%
\end{equation}

\end{definition}

For the strong moderation, more precise estimates are given for the constants
which appear in (\ref{Defmod1}).

\begin{proposition}
\label{PropExt}~\newline(i)~Any moderate net $\left(  L_{\varepsilon}\right)
_{\varepsilon}\in\left(  \mathcal{L}\left(  \mathcal{S}\left(  \mathbb{R}%
^{n}\right)  ,\mathcal{O}_{M}\left(  \mathbb{R}^{m}\right)  \right)  \right)
^{\left(  0,1\right]  }$ can be extended to a map $L$ belonging to
$\mathbf{L(}\mathcal{G}_{\mathcal{S}}\left(  \mathbb{R}^{n}\right)
,\mathcal{G}_{\tau}\left(  \mathbb{R}^{m}\right)  )$ and defined by
\begin{equation}
L\left(  f\right)  =\left(  L_{\varepsilon}\left(  f_{\varepsilon}\right)
\right)  _{\varepsilon}+\mathcal{N}_{\tau}\left(  \mathbb{R}^{m}\right)
,\label{GSTextDef}%
\end{equation}
where $\left(  f_{\varepsilon}\right)  _{\varepsilon}$ is any representative
of $f$.\newline(ii)~The extension $L$ depends only on the family $\left(
L_{\varepsilon}\right)  _{\varepsilon}$ in the following sense: If $\left(
N_{\varepsilon}\right)  _{\varepsilon}$ is a negligible net of maps, then the
extensions of $\left(  L_{\varepsilon}\right)  _{\varepsilon}$ and $\left(
L_{\varepsilon}+N_{\varepsilon}\right)  _{\varepsilon}$ are equal.\newline%
(iii)~If the family $\left(  L_{\varepsilon}\right)  _{\varepsilon}$ is
moderate, with the assumption that the net of constants $\left(
C_{\varepsilon}\right)  _{\varepsilon}$ in (\ref{Defmod1}) satisfies
$C_{\varepsilon}=\mathrm{O}\left(  \varepsilon^{-r(l)}\right)  $ with
$\underset{l\rightarrow+\infty}{\lim\sup}\left(  r(l)/l\right)  <c$, then
$L\left(  \mathcal{G}_{\mathcal{S}}^{\infty}\left(  \mathbb{R}^{n}\right)
\right)  $ is included in $\mathcal{G}_{\tau}^{\mathcal{L}_{c}}\left(
\mathbb{R}^{m}\right)  $.\newline(iv)~Let $\left(  a,b,c\right)  $ be in
$\left(  \mathbb{R}_{+}\right)  ^{3}$: If the net $\left(  L_{\varepsilon
}\right)  _{\varepsilon}$ is $\mathcal{L}_{b,c}$-strongly moderate, then
$L\left(  \mathcal{G}_{\mathcal{S}}^{\mathcal{L}_{a}}\left(  \mathbb{R}%
^{n}\right)  \right)  $ is included in $\mathcal{G}_{\tau}^{\mathcal{L}%
_{ab+c}}\left(  \mathbb{R}^{m}\right)  $.\newline Moreover, $L\left(
\mathcal{G}_{\mathcal{S}}^{\mathcal{L}_{0}}\left(  \mathbb{R}^{n}\right)
\right)  $ is included in $\mathcal{G}_{\tau}^{\mathcal{L}_{c}}\left(
\mathbb{R}^{m}\right)  $ even if $b=+\infty$.
\end{proposition}

\begin{proof}
\textit{Assertions (i)} \& \textit{(ii).}-~Fix $l\in\mathbb{N}$ and let
$\left(  f_{\varepsilon}\right)  _{\varepsilon}$ be in $\mathcal{E}%
_{\mathcal{S}}\left(  \mathbb{R}^{n}\right)  $. According to the definition of
moderate nets, we get $\left(  C_{\varepsilon}\right)  _{\varepsilon}%
\in\mathcal{E}_{M}\left(  \mathbb{R}_{+}\right)  $ and $\left(  p,q,l^{\prime
}\right)  \in\mathbb{N}^{3}$ such that
\begin{equation}
\mu_{-p,l}\left(  L_{\varepsilon}\left(  f_{\varepsilon}\right)  \right)  \leq
C_{\varepsilon}\,\mu_{q,l^{\prime}}\left(  f_{\varepsilon}\right)  \text{, for
}\varepsilon\text{ small enough}. \label{GSTPrExt1}%
\end{equation}
Inequality (\ref{GSTPrExt1}) leads to $\left(  L_{\varepsilon}\left(
f_{\varepsilon}\right)  \right)  _{\varepsilon}\in\mathcal{E}_{\tau}\left(
\mathbb{R}^{m}\right)  $. Moreover, if $\left(  f_{\varepsilon}\right)
_{\varepsilon}$ belongs to $\mathcal{N}_{\mathcal{S}}\left(  \mathbb{R}%
^{n}\right)  $, the same inequality implies that $\left(  L_{\varepsilon
}\left(  f_{\varepsilon}\right)  \right)  _{\varepsilon}\in\mathcal{N}_{\tau
}\left(  \mathbb{R}^{m}\right)  $. These two properties show that $L$ is well
defined by formula (\ref{GSTextDef}). The $\overline{\mathbb{C}}$ linearity
follows from from the moderation. Moreover, inequality (\ref{GSTPrExt1})
implies easily the second assertion.\smallskip

\noindent\textit{Assertion (iii)} \& \textit{(iv).}-~We shall prove
\textit{(iv)} for $a\in\left(  0,+\infty\right)  $,\ since the proof of
\textit{(iii)} and of the case $a=0$ in \textit{(iv) }are even simpler.
Suppose that $\left(  L_{\varepsilon}\right)  _{\varepsilon}$ is
$\mathcal{L}_{b,c}$-strongly moderate and consider $\left(  f_{\varepsilon
}\right)  _{\varepsilon}\in\mathcal{E}_{\mathcal{S}}^{\mathcal{L}_{a}}\left(
\mathbb{R}^{n}\right)  $. There exist a sequence $\lambda\in\mathbb{R}%
_{+}^{\mathbb{N}}$, with $\lim\!\sup_{l\rightarrow+\infty}\left(
\lambda(l)/l\right)  <b$, and a sequence $r\in\mathbb{R}_{+}^{\mathbb{N}}$,
with $\lim\!\sup_{l\rightarrow+\infty}\left(  r(l)/l\right)  <c$, such that
\[
\forall l\in\mathbb{N},\;\exists C\in\mathbb{R}_{+},\ \exists\left(
p,q\right)  \in\mathbb{N}^{2},\;\;\mu_{-p,l}\left(  L_{\varepsilon}\left(
f_{\varepsilon}\right)  \right)  \leq C\,\varepsilon^{-r(l)}\mu_{q,\lambda
(l)}\left(  f_{\varepsilon}\right)  \text{ (for }\varepsilon\text{ small
enough)}.
\]
As $\left(  f_{\varepsilon}\right)  _{\varepsilon}$ is in $\mathcal{E}%
_{\mathcal{S}}^{\mathcal{L}_{a}}\left(  \mathbb{R}^{n}\right)  $, there exists
a sequence $N\in\mathbb{R}_{+}^{\mathbb{N}},\;$with $\lim\!\sup_{\lambda
\rightarrow+\infty}\left(  N(\lambda)/\lambda\right)  <a$, such that
\[
\forall s\in\mathbb{N},\ \forall\lambda\in\mathbb{N},\;\;\mu_{s,\lambda
}\left(  f_{\varepsilon}\right)  =\mathrm{O}\left(  \varepsilon^{-N(\lambda
)}\right)  \;\mathrm{as}\;\varepsilon\rightarrow0.
\]
We get that
\[
\forall l\in\mathbb{N},\;\exists C\in\mathbb{R}_{+},\ \exists p\in
\mathbb{N},\;\;\mu_{-p,l}\left(  L_{\varepsilon}\left(  f_{\varepsilon
}\right)  \right)  \leq C\,\varepsilon^{-N_{1}\left(  l\right)  }%
\;\;\text{\ with }N_{1}\left(  l\right)  =r(l)+N\left(  \lambda(l)\right)  ,
\]
for $\varepsilon$ small enough.

\begin{itemize}
\item If $\lambda(l)$ is bounded, we immediately have: $N_{1}\left(  l\right)
/l=\mathrm{O}\left(  r(l)/l\right)  $ for $l\rightarrow+\infty$.

\item If $\lambda(l)$ is not bounded, for $\lambda(l)\neq0$, we have%
\begin{equation}
\frac{N_{1}\left(  l\right)  }{l}=\frac{r(l)}{l}+\frac{N\left(  \lambda
(l)\right)  }{\lambda(l)}\frac{\lambda(l)}{l}. \label{GSTmodmapsXX}%
\end{equation}
We have $\lim\sup_{l\rightarrow+\infty}\left(  N\left(  \lambda(l)\right)
/\lambda(l)\right)  <a$ and thus $\lim\sup_{l\rightarrow+\infty}\frac{N\left(
\lambda(l)\right)  }{\lambda(l)}\frac{\lambda(l)}{l}<ab$. This gives
\[
\underset{l\rightarrow+\infty}{\lim\sup}\left(  N_{1}(l)/l\right)  <ab+c
\]
and $\left(  L_{\varepsilon}\left(  f_{\varepsilon}\right)  \right)
_{\varepsilon}\in\mathcal{E}_{\tau}^{\mathcal{L}_{ab+c}}\left(  \mathrm{C}%
^{\infty}\left(  \mathbb{R}^{m}\right)  \right)  $, which shows the assertion.
\end{itemize}

Finally, if $\left(  f_{\varepsilon}\right)  _{\varepsilon}$ is in
$\mathcal{E}_{\mathcal{S}}^{\mathcal{L}_{0}}\left(  \mathbb{R}^{n}\right)  $,
the sequence $N$ can be chosen such that $\lim_{\lambda\rightarrow+\infty
}\left(  N\left(  \lambda\right)  /\lambda\right)  =0$. Then, for $b=+\infty$,
the sequence $l\mapsto\lambda(l)/l$ is bounded. It follows that $\lim
\!\sup_{l\rightarrow+\infty}\left(  N_{1}(l)/l\right)  <c$.
\end{proof}

\subsection{Main results}

\begin{theorem}
\label{KernSTHKern1}Consider $\left(  a,b,c\right)  \in\mathbb{R}_{+}^{3}$
with $a\leq1$ and $ab+c\leq1$. Let $\left(  L_{\varepsilon}\right)
_{\varepsilon}\in\mathcal{L}\left(  \mathcal{S}(\mathbb{R}^{n}),\mathcal{O}%
_{M}(\mathbb{R}^{m})\right)  ^{\left(  0,1\right]  }$ be a net of
$\mathcal{L}_{b,c}$-strongly moderate linear maps and $L\in\mathbf{L}\left(
\mathcal{G}_{\mathcal{S}}(\mathbb{R}^{n}),\mathcal{G}_{\tau}(\mathbb{R}%
^{m})\right)  $ its canonical extension. There exists $H_{L}\in\mathcal{G}%
_{\tau}\left(  \mathbb{R}^{m+n}\right)  $ such that
\begin{equation}
\forall f\in\mathcal{G}_{\mathcal{S}}^{\mathcal{L}_{_{a}}}\left(
\mathbb{R}^{n}\right)  ,\;\;L\left(  f\right)  =\left[  \left(  x\longmapsto
\int H_{L,\varepsilon}(x,y)f_{\varepsilon}(y)\,\mathrm{d}y\right)
_{\varepsilon}\right]  ,\label{GSTScequa1}%
\end{equation}
where $\left(  H_{L,\varepsilon}\right)  _{\varepsilon}$ (\textit{resp.}
$\left(  f_{\varepsilon}\right)  _{\varepsilon}$) is any representative of
$H_{L}$ (\textit{resp.} $f$).
\end{theorem}

\begin{remark}
In Theorem \ref{KernSTHKern1}, the parameter $b$ (resp. $c$ ) is related to
the \textquotedblleft regularity\textquotedblright\ of the net $\left(
L_{\varepsilon}\right)  _{\varepsilon}$, with respect to the derivative index
$l$ in the family of semi-norms $\left(  \mu_{r,l}\right)  _{r,l}$ (resp. to
the parameter $\varepsilon$). The more \textquotedblleft
irregular\textquotedblright\ the net of maps $\left(  L_{\varepsilon}\right)
_{\varepsilon}$ is, that is the bigger $b$ is (resp. the closer to $1$ $c$
is), the smaller is the space on which equality (\ref{GSTScequa1}) holds. The
limit cases for $c$ are $c=1$ (for which $a=0$ and (\ref{GSTScequa1}) holds
only on $\mathcal{G}_{\mathcal{S}}^{\mathcal{L}_{_{0}}}\left(  \mathbb{R}%
^{n}\right)  $) and $c=0$ (the net of constants $\left(  C_{\varepsilon
}\right)  _{\varepsilon}$ in relation (\ref{Defmod1}) depends slowly on
$\varepsilon$) for which the conditions on $\left(  a,b,c\right)  $ are
reduced to $a<1$ and $ab\leq1$. (Note that these limiting conditions are
induced by Lemma \ref{LmnMTh4}.) \medskip
\end{remark}

By using Proposition \ref{PropExt}-\emph{(iii)}, we can obtain an analogon of
Theorem \ref{KernSTHKern1} valid for more irregular nets of maps.

\begin{theorem}
\label{KernSTHKern12}Let $\left(  L_{\varepsilon}\right)  _{\varepsilon}%
\in\mathcal{L}\left(  \mathcal{S}(\mathbb{R}^{n}),\mathcal{O}_{M}%
(\mathbb{R}^{m})\right)  ^{\left(  0,1\right]  }$ be a net of moderate linear
maps such that the net of constants $\left(  C_{\varepsilon}\right)
_{\varepsilon}$ in relation (\ref{Defmod1}) satisfies $C_{\varepsilon
}=\mathrm{O}\left(  \varepsilon^{-r(l)}\right)  $ with $r\in\mathcal{L}_{_{1}%
}$. Then, the extension $\left(  L_{\varepsilon}\right)  _{\varepsilon}$
admits an integral representation such that relation (\ref{GSTScequa1}) holds
for $f$ in $\mathcal{G}_{\mathcal{S}}^{\infty}\left(  \mathbb{R}^{n}\right)  $.
\end{theorem}

We turn now to the relationship with the classical isomorphism result:
Consider
\[
\Lambda\in\mathcal{L}\left(  \mathcal{S}(\mathbb{R}^{n}),\mathcal{S}^{\prime
}(\mathbb{R}^{m})\right)
\]
and define a net of linear mappings $\left(  L_{\varepsilon}\right)
_{\varepsilon}$ by
\[
L_{\varepsilon}:\ \ \mathcal{S}\left(  \mathbb{R}^{n}\right)  \rightarrow
\mathrm{C}^{\infty}\left(  \mathbb{R}^{m}\right)  ,\;\;f\mapsto\Lambda\left(
f\right)  \ast\varphi_{\varepsilon^{s}}\text{,\ \ (}s\text{ fixed real
parameter in }\left(  0,1\right)  \text{)}%
\]
where $\left(  \varphi_{\varepsilon}\right)  _{\varepsilon}$ is defined as in
(\ref{CSTGSBinfini2}), starting from $\varphi\in\mathcal{S}\left(
\mathbb{R}^{m}\right)  $ which satisfies (\ref{CSTGSBinfini}). We have:

\begin{proposition}
\label{GSTCoroTHS}~\newline(i)~For all $\varepsilon\in\left(  0,1\right]  $,
$L_{\varepsilon}$ is continuous for the usual topologies of $\mathcal{S}%
\left(  \mathbb{R}^{n}\right)  $ and $\mathcal{O}_{M}\left(  \mathbb{R}%
^{m}\right)  $ and the net $\left(  L_{\varepsilon}\right)  _{\varepsilon}$ is
$\left(  0,s\right)  $-strongly moderate.\newline(ii)~From (i), the extension
$L$ of the net $\left(  L_{\varepsilon}\right)  _{\varepsilon}$ admits a
kernel $H_{L}$.\ Furthermore, for all $f\in\mathcal{S}\left(  \mathbb{R}%
^{n}\right)  $, $\Lambda\left(  f\right)  $ is equal to $\widetilde{H}%
_{L}\left(  f\right)  $ in the \emph{generalized distribution sense
}\cite{NePiSc}, that is
\begin{equation}
\forall g\in\mathcal{S}\left(  \mathbb{R}^{m}\right)  ,\;\;\;\left\langle
\Lambda\left(  f\right)  ,g\right\rangle =\left\langle \widetilde{H}%
_{L}\left(  f\right)  ,g\right\rangle \text{ in }\overline{\mathbb{C}}\text{.}
\label{KernsEquaGTD}%
\end{equation}

\end{proposition}

In other words, equality (\ref{KernsEquaGTD}) means that, for all
$p\in\mathbb{N}$,
\begin{equation}
\forall g\in\mathcal{S}\left(  \mathbb{R}^{m}\right)  ,\;\;\left\langle
\Lambda\left(  f\right)  ,g\right\rangle -\int\left(  \int H_{L,\varepsilon
}\left(  x,y\right)  f\left(  y\right)  \,\mathrm{d}y\right)  g\left(
x\right)  \,\mathrm{d}x=\mathrm{O}\left(  \varepsilon^{p}\right)  \text{, as
}\varepsilon\rightarrow0, \label{GSTDefgdequa}%
\end{equation}
where $\left(  H_{L,\varepsilon}\right)  _{\varepsilon}$ is any representative
of $H_{L}$. In particular, this result implies that $\Lambda\left(  f\right)
$ and $\widetilde{H}_{L}\left(  f\right)  $ are associated or weakly equal,
\emph{i.e.}%
\[
\left(  x\mapsto\int H_{L,\varepsilon}\left(  x,y\right)  f\left(  y\right)
\,\mathrm{d}y\right)  \longrightarrow\Lambda\left(  f\right)  \text{ in
}\mathcal{S}^{\prime}\left(  \mathbb{R}^{m}\right)  \text{ as }\varepsilon
\rightarrow0\text{.}%
\]

\section{Proofs of Theorem \ref{KernSTHKern1} and Proposition \ref{GSTCoroTHS}%
}

\subsection{Proof of Theorem \ref{KernSTHKern1}}

We shall prove Theorem \ref{KernSTHKern1}. (The proof of Theorem
\ref{KernSTHKern12} follows the same lines.) Consider $\varphi\in
\mathcal{S}\left(  \mathbb{R}^{m}\right)  $ (resp. $\psi\in\mathcal{S}\left(
\mathbb{R}^{n}\right)  $) which satisfies (\ref{CSTGSBinfini}) and define
$\left(  \varphi_{\varepsilon}\right)  _{\varepsilon}$ (resp. $\left(
\psi_{\varepsilon}\right)  _{\varepsilon}$) as in (\ref{CSTGSBinfini2}).\ For
all $\varepsilon\in\left(  0,1\right]  $ and $y\in\mathbb{R}^{n}$, we set
\[
\psi_{\varepsilon,\,\cdot}:\ \ \mathbb{R}^{n}\rightarrow\mathcal{S}\left(
\mathbb{R}^{n}\right)  ,\;\;y\mapsto\psi_{\varepsilon,y}=\left\{  v\mapsto
\psi_{\varepsilon}\left(  y-v\right)  \right\}  .
\]
Then, for all $\varepsilon\in\left(  0,1\right]  $, the map
\[
H_{\varepsilon}:\ \ \mathbb{R}^{m+n}\rightarrow\mathbb{C},\ \ \left(
x,y\right)  \mapsto\left(  L_{\varepsilon}\left(  \psi_{\varepsilon,y}\right)
\ast\varphi_{\varepsilon}\right)  \left(  x\right)  =\int L_{\varepsilon
}\left(  \psi_{\varepsilon,y}\right)  \left(  x-u\right)  \varphi
_{\varepsilon}\left(  u\right)  \,\mathrm{d}u,
\]
is well defined.

Indeed, $L_{\varepsilon}\left(  \psi_{\varepsilon,y}\right)  $ belongs to
$\mathcal{O}_{M}\left(  \mathbb{R}\right)  $ and $\varphi_{\varepsilon}$ to
$\mathcal{S}\left(  \mathbb{R}^{m}\right)  $, making $L_{\varepsilon}\left(
\psi_{\varepsilon,y}\right)  \left(  x-\cdot\right)  \varphi_{\varepsilon
}\left(  \cdot\right)  $ -\thinspace and its derivatives\thinspace-
\textrm{L}$^{1}$ functions.

\begin{lemma}
\label{LmnMTh1}For all $\varepsilon\in\left(  0,1\right]  $, $H_{\varepsilon}$
is of class $\mathrm{C}^{\infty}$ and $\left(  H_{\varepsilon}\right)
_{\varepsilon}\in\mathcal{E}_{\tau}\left(  \mathbb{R}^{m+n}\right)  $.
\end{lemma}

\begin{proof}
The fact that $H_{\varepsilon}$ is of class $\mathrm{C}^{\infty}$ follows from
classical arguments of integral calculus.\ It also uses the topological
isomorphism between $\mathrm{C}^{\infty}\left(  \mathbb{R}^{d_{1}}%
,\mathrm{C}^{\infty}\left(  \mathbb{R}^{d_{2}}\right)  \right)  $ and
$\mathrm{C}^{\infty}\left(  \mathbb{R}^{d_{1}+d_{2}},\mathbb{C}\right)  $
($d_{1},$ $d_{2}\in\mathbb{N}\backslash\left\{  0\right\}  $), the linearity
and continuity of both $L_{\varepsilon}$ and the convolution. (See Lemma 28
and 29 in \cite{ADSKTT} for very close proofs.)

Let us now consider $\left(  \alpha,\beta\right)  \in\mathbb{N}^{m+n}$ and
$\partial_{x}^{\alpha}$ (\textit{resp}. $\partial_{y}^{\beta}$) the $\alpha
$-partial derivative (\textit{resp}. $\beta$-partial derivative) with respect
to the variable $x$ (\textit{resp}. $y$). and set $l=\left\vert \beta
\right\vert $. We have
\begin{align}
\left.  \forall\left(  x,y\right)  \in\mathbb{R}^{m+n},\ \ \ \partial
_{x}^{\alpha}\partial_{y}^{\beta}H_{\varepsilon}\left(  x,y\right)  \right.
&  =\left(  \partial_{y}^{\beta}L_{\varepsilon}\left(  \psi_{\varepsilon
,y}\right)  \ast\partial_{x}^{\alpha}\varphi_{\varepsilon}\right)  \left(
x\right) \nonumber\\
&  =\int\partial_{y}^{\beta}L_{\varepsilon}\left(  \psi_{\varepsilon
,y}\right)  \left(  x-u\right)  \partial_{x}^{\alpha}\varphi_{\varepsilon
}\left(  u\right)  \,\mathrm{d}u\nonumber\\
&  =\int\partial_{y}^{\beta}L_{\varepsilon}\left(  \psi_{\varepsilon
,y}\right)  \left(  x-\varepsilon\zeta\right)  \partial_{x}^{\alpha}%
\varphi\left(  \zeta\right)  \,\mathrm{d}\zeta\label{KernSHmod1}%
\end{align}
Using the moderation of $\left(  L_{\varepsilon}\right)  _{\varepsilon}$, we
get the existence of $\left(  C_{\varepsilon}\right)  _{\varepsilon}%
\in\mathcal{E}_{M}\left(  \mathbb{R}_{+}\right)  $,\ \ $\left(  p,q,l^{\prime
}\right)  \in\mathbb{N}^{3}$ such that, for $\varepsilon$ small enough,%
\begin{align*}
\left.  \forall\left(  x,\zeta\right)  \in\mathbb{R}^{2m},\ \ \ \left\vert
\partial_{y}^{\beta}L_{\varepsilon}\left(  \psi_{\varepsilon,y}\right)
\left(  x-u\right)  \right\vert \right.   &  \leq C_{\varepsilon}\left(
1+\left\vert x-\varepsilon\zeta\right\vert \right)  ^{p}\mu_{q,l^{\prime}%
}\left(  \psi_{\varepsilon,y}\right) \\
&  \leq C_{\varepsilon}\left(  1+\left\vert x\right\vert \right)  ^{p}\left(
1+\left\vert \zeta\right\vert \right)  ^{p}\mu_{q,l^{\prime}}\left(
\psi_{\varepsilon,y}\right)  .
\end{align*}
We have
\begin{multline*}
\mu_{q,l^{\prime}}\left(  \psi_{\varepsilon,y}\right)  =\sup_{w\in
\mathbb{R}^{n},\ \left\vert \alpha\right\vert \leq l^{\prime}}\left(
1+\left\vert y-w\right\vert \right)  ^{q}\left\vert \partial^{\alpha}%
\psi_{\varepsilon}(w)\right\vert \leq\left(  1+\left\vert y\right\vert
\right)  ^{q}\mu_{q,l^{\prime}}\left(  \psi_{\varepsilon}\right) \\
\leq\varepsilon^{-n-l^{\prime}}\left(  1+\left\vert y\right\vert \right)
^{q}\mu_{q,l^{\prime}}\left(  \psi\right)  .
\end{multline*}
Therefore, there exists $\left(  C_{\varepsilon}^{\prime}\right)
_{\varepsilon}\in\mathcal{E}_{M}\left(  \mathbb{R}_{+}\right)  $ such that
\[
\left.  \forall\left(  x,\zeta,y\right)  \in\mathbb{R}^{2m+n},\ \ \ \left\vert
\partial_{y}^{\beta}L_{\varepsilon}\left(  \psi_{\varepsilon,y}\right)
\left(  x-\varepsilon\zeta\right)  \right\vert \right.  \leq C_{\varepsilon
}^{\prime}\left(  1+\left\vert x\right\vert \right)  ^{p}\left(  1+\left\vert
\zeta\right\vert \right)  ^{p}\left(  1+\left\vert y\right\vert \right)
^{q}.
\]
for $\varepsilon$ small enough. As $\partial_{x}^{\alpha}\varphi$ is rapidly
decreasing, by replacing the last estimate above in (\ref{KernSHmod1}), we get
the existence of $C_{\varepsilon}^{\prime\prime}\in\mathcal{E}_{M}\left(
\mathbb{R}_{+}\right)  $ such that
\[
\forall\left(  x,y\right)  \in\mathbb{R}^{m+n},\ \ \ \left\vert \partial
_{x}^{\alpha}\partial_{y}^{\beta}H_{\varepsilon}\left(  x,y\right)
\right\vert \leq C_{\varepsilon}^{\prime\prime}\left(  1+\left\vert
x\right\vert \right)  ^{p}\left(  1+\left\vert y\right\vert \right)  ^{q}\leq
C_{\varepsilon}^{\prime\prime}\left(  1+\left\vert \left(  x,y\right)
\right\vert \right)  ^{p+q},
\]
for $\varepsilon$ small enough. Thus, $\left(  H_{\varepsilon}\right)
_{\varepsilon}\in\mathcal{E}_{\tau}\left(  \mathbb{R}^{m+n}\right)  $ as claimed.
\end{proof}

\begin{lemma}
\label{LmnMTh3}For all $\left(  f_{\varepsilon}\right)  _{\varepsilon}$ in
$\mathcal{E}_{S}\left(  \mathbb{R}^{n}\right)  $, we have
\[
\widetilde{H}_{\varepsilon}\left(  f_{\varepsilon}\right)  \left(  x\right)
=\left(  L_{\varepsilon}\left(  \psi_{\varepsilon}\ast f_{\varepsilon}\right)
\ast\varphi_{\varepsilon}\right)  (x).
\]

\end{lemma}

\begin{proof}
Let $\left(  f_{\varepsilon}\right)  _{\varepsilon}$ be in $\mathcal{E}%
_{S}\left(  \mathbb{R}^{n}\right)  $. For any fixed $\varepsilon\in\left(
0,1\right]  $ and $x\in\mathbb{R}^{m}$, we have%
\[
\widetilde{H}_{\varepsilon}\left(  f_{\varepsilon}\right)  \left(  x\right)
=\int\left(  \int L_{\varepsilon}\left(  \psi_{\varepsilon,y}\right)  \left(
x-u\right)  \varphi_{\varepsilon}\left(  u\right)  \,\mathrm{d}u\right)
f_{\varepsilon}(y)\,\mathrm{d}y.
\]
Using a similar argument as in the proof of Lemma \ref{LmnMTh1}, we get the
existence of $C_{\varepsilon}(x)>0$ such that
\[
\forall u\in\mathbb{R}^{m},\ \ \left\vert L_{\varepsilon}\left(
\psi_{\varepsilon,y}\right)  \left(  x-u\right)  \right\vert \leq
C_{\varepsilon}(x)\left(  1+\left\vert u\right\vert \right)  ^{p}.
\]
Thus, the map $\left(  u,y\right)  \mapsto L_{\varepsilon}\left(
\psi_{\varepsilon,y}\right)  \left(  x-u\right)  \varphi_{\varepsilon}\left(
u\right)  f_{\varepsilon}(y)$ is in \textrm{L}$^{1}\left(  \mathbb{R}%
^{m+n}\right)  $ and, by\ Fubini's Theorem,%
\begin{align*}
\widetilde{H}_{\varepsilon}\left(  f_{\varepsilon}\right)  \left(  x\right)
&  =\int\left(  \int L_{\varepsilon}\left(  \psi_{\varepsilon,y}\right)
\left(  x-u\right)  f_{\varepsilon}(y)\,\mathrm{d}y\right)  \varphi
_{\varepsilon}\left(  u\right)  \,\mathrm{d}u\\
&  =(\,\{\xi\mapsto\int L_{\varepsilon}\left(  \psi_{\varepsilon,y}\right)
\left(  \xi\right)  f_{\varepsilon}(y)\,\mathrm{d}y\}\ast\varphi_{\varepsilon
}\,)\left(  x\right)  .
\end{align*}
An adaptation of the proof of Lemma 30 in \cite{ADSKTT} shows that, for all
$\xi\in\mathbb{R}^{m}$, we have the following equality
\begin{align}
\left.  \forall g\in\mathcal{D}\left(  \mathbb{R}^{n}\right)  ,\ \ \ \int
L_{\varepsilon}\left(  \psi_{\varepsilon,y}\right)  \left(  \xi\right)
g(y)\,\mathrm{d}y\right.   &  =L_{\varepsilon}\left(  \{v\mapsto\int
\psi_{\varepsilon,y}\left(  v\right)  g(y)\mathrm{d}y\}\right)  \left(
\xi\right) \label{KernSInter1}\\
(  &  =L_{\varepsilon}\left(  \{v\mapsto\int\psi_{\varepsilon}\left(
y-v\right)  g(y)\mathrm{d}y\}\right)  \left(  \xi\right)  \,).\nonumber
\end{align}
(Indeed, the integrals under consideration in (\ref{KernSInter1}) are
integrals of continuous functions on compact sets and can be considered as
limits of Riemann sums in the spirit of \cite{HorPDOT1}, Lemma 4.1.3. The
linearity and the continuity of $L_{\varepsilon}$ allows to exchange the order
or the operations integral and $L_{\varepsilon}$.)

Then, a density argument shows that equality (\ref{KernSInter1}) holds for
$g\in\mathcal{S}\left(  \mathbb{R}^{n}\right)  $.\ Thus%
\[
\int L_{\varepsilon}\left(  \psi_{\varepsilon,y}\right)  \left(  \xi\right)
f_{\varepsilon}(y)\,\mathrm{d}y=L_{\varepsilon}\left(  \{v\mapsto\int
\psi_{\varepsilon}\left(  y-v\right)  f_{\varepsilon}(y)\mathrm{d}y\}\right)
\left(  \xi\right)  =L_{\varepsilon}\left(  \psi_{\varepsilon}\ast
f_{\varepsilon}\right)  \left(  \xi\right)
\]
and $\widetilde{H}_{\varepsilon}\left(  f_{\varepsilon}\right)  \left(
x\right)  =\left(  L_{\varepsilon}\left(  \psi_{\varepsilon}\ast
f_{\varepsilon}\right)  \ast\varphi_{\varepsilon}\right)  (x)$ as
claimed.\medskip\smallskip
\end{proof}

We now complete the proof of Theorem \ref{KernSTHKern1}.\ Set
\[
H_{L}=\left[  \left(  H_{\varepsilon}\right)  _{\varepsilon}\right]  _{\tau
}=\left(  \left(  x,y\right)  \mapsto\left(  \,\Psi_{\varepsilon,y}\ast
\varphi_{\varepsilon}\right)  \left(  x\right)  \,\right)  _{\varepsilon
}+\mathcal{N}_{\tau}\left(  \mathbb{R}^{m+n}\right)  .
\]
For all $\left(  f_{\varepsilon}\right)  _{\varepsilon}$ in $\mathcal{E}%
_{\mathcal{S}}^{\mathcal{L}_{a}}\left(  \mathbb{R}^{n}\right)  $, we have
\[
\widetilde{H}_{L}\left(  \left[  \left(  f_{\varepsilon}\right)
_{\varepsilon}\right]  _{\mathcal{S}}\right)  =\left[  \left(  \widetilde
{H}_{\varepsilon}\left(  f_{\varepsilon}\right)  \right)  _{\varepsilon
}\right]  _{\tau},
\]
by definition of the integral operator.\ We have to compare $\left(
\widetilde{H}_{\varepsilon}\left(  f_{\varepsilon}\right)  \right)
_{\varepsilon}$ and $\left(  L_{\varepsilon}\left(  f_{\varepsilon}\right)
\right)  _{\varepsilon}$. According to Lemma \ref{LmnMTh3}, we have for all
$\varepsilon\in\left(  0,1\right]  $
\begin{align*}
\widetilde{H}_{\varepsilon}\left(  f_{\varepsilon}\right)  -L_{\varepsilon
}\left(  f_{\varepsilon}\right)   &  =\left(  L_{\varepsilon}\left(
\psi_{\varepsilon}\ast f_{\varepsilon}\right)  \ast\varphi_{\varepsilon
}\right)  -L_{\varepsilon}\left(  f_{\varepsilon}\right) \\
&  =L_{\varepsilon}\left(  \psi_{\varepsilon}\ast f_{\varepsilon}\right)
\ast\varphi_{\varepsilon}-L_{\varepsilon}\left(  f_{\varepsilon}\right)
\ast\varphi_{\varepsilon}+L_{\varepsilon}\left(  f_{\varepsilon}\right)
\ast\varphi_{\varepsilon}-L_{\varepsilon}\left(  f_{\varepsilon}\right) \\
&  =L_{\varepsilon}\left(  \psi_{\varepsilon}\ast f_{\varepsilon
}-f_{\varepsilon}\right)  \ast\varphi_{\varepsilon}+L_{\varepsilon}\left(
f_{\varepsilon}\right)  \ast\varphi_{\varepsilon}-L_{\varepsilon}\left(
f_{\varepsilon}\right)  .
\end{align*}
Remarking that $\left(  f_{\varepsilon}\right)  _{\varepsilon}\in
\mathcal{E}_{\mathcal{S}}^{\mathcal{L}_{a}}\left(  \mathbb{R}^{n}\right)  $
and $\left(  L_{\varepsilon}\left(  f_{\varepsilon}\right)  \right)
_{\varepsilon}\in\mathcal{E}_{\tau}^{\mathcal{L}_{a+bc}}\left(  \mathbb{R}%
^{m}\right)  \subset\mathcal{E}_{\tau}^{\mathcal{L}_{1}}\left(  \mathbb{R}%
^{m}\right)  $, we get
\[
\left(  \,L_{\varepsilon}\left(  f_{\varepsilon}\right)  \ast\varphi
_{\varepsilon}-L_{\varepsilon}\left(  f_{\varepsilon}\right)  \,\right)
_{\varepsilon}\in\mathcal{N}_{\tau}\left(  \mathbb{R}^{m}\right)
\;\;\mathrm{and}\;\;\left(  \,\psi_{\varepsilon}\ast f_{\varepsilon
}-f_{\varepsilon}\,\right)  _{\varepsilon}\in\mathcal{N}_{\emph{S}}\left(
\mathbb{R}^{m}\right)
\]
by Lemma \ref{LmnMTh4}. This last property gives
\[
\left(  \,L_{\varepsilon}\left(  \psi_{\varepsilon}\ast f_{\varepsilon
}-f_{\varepsilon}\right)  \,\right)  _{\varepsilon}\in\mathcal{N}_{\tau
}\left(  \mathbb{R}^{m}\right)  \;\;\mathrm{and}\;\;\left(  \,L_{\varepsilon
}\left(  \psi_{\varepsilon}\ast f_{\varepsilon}-f_{\varepsilon}\right)
\ast\varphi_{\varepsilon}\,\right)  \,\in\mathcal{N}_{\tau}\left(
\mathbb{R}^{m}\right)  ,
\]
since $\left(  \eta_{\varepsilon}\ast\varphi_{\varepsilon}\right)
_{\varepsilon}\in\mathcal{N}_{\tau}\left(  \mathbb{R}^{m}\right)  $ for all
$\left(  \eta_{\varepsilon}\right)  _{\varepsilon}\in\mathcal{N}_{\tau}\left(
\mathbb{R}^{m}\right)  $. Finally
\[
\left[  \left(  \widetilde{H}_{\varepsilon}\left(  f_{\varepsilon}\right)
\right)  _{\varepsilon}\right]  _{\tau}=\left[  \left(  L_{\varepsilon}\left(
f_{\varepsilon}\right)  \right)  _{\varepsilon}\right]  _{\tau}=L\left(
\left[  \left(  f_{\varepsilon}\right)  _{\varepsilon}\right]  _{\mathcal{S}%
}\right)  ,
\]
this last equality by definition of the extension of a linear map.

\subsection{Proof of Proposition \ref{GSTCoroTHS}}

\textit{Assertion }$\left(  i\right)  $.-~For a fixed $\varepsilon\in\left(
0,1\right]  $, $L_{\varepsilon}$ is obtained by composition of the continuous
maps $\Lambda:\mathcal{S}\left(  \mathbb{R}^{n}\right)  \mapsto\mathcal{S}%
^{\prime}\left(  \mathbb{R}^{m}\right)  $ and
\[
\mathcal{S}^{\prime}\left(  \mathbb{R}^{m}\right)  \rightarrow\mathcal{O}%
_{M}\left(  \mathbb{R}^{n}\right)  ,\ \ T\mapsto T_{k}\ast\varphi
_{\varepsilon^{s}}%
\]
Thus $L_{\varepsilon}$ is continuous. We have now to show that the net
$\left(  L_{\varepsilon}\right)  _{\varepsilon}\in\mathcal{L}\left(
\mathcal{S}\left(  \mathbb{R}^{n}\right)  ,\mathcal{O}_{M}\left(
\mathbb{R}^{m}\right)  \right)  ^{\left(  0,1\right]  }$ is strongly moderate.
We have
\begin{align*}
\left.  \forall f\in\mathcal{S}\left(  \mathbb{R}^{n}\right)  ,\;\forall
x\in\mathbb{R}^{m},\;\forall\alpha\in\mathbb{N}^{m},\;\;\;\partial^{\alpha
}\left(  L_{\varepsilon}(f)\right)  (x)\right.   &  =\left(  \Lambda\left(
f\right)  \ast\partial^{\alpha}\varphi_{\varepsilon^{s}}\right)  (x)\\
&  =\left\langle \Lambda\left(  f\right)  ,\left\{  y\mapsto\partial^{\alpha
}\varphi_{\varepsilon^{s}}\left(  x-y\right)  \right\}  \right\rangle .
\end{align*}
The map
\[
\Theta:\mathcal{S}\left(  \mathbb{R}^{n}\right)  \times\mathcal{S}\left(
\mathbb{R}^{m}\right)  ,\;\;\;\left(  f,\varphi\right)  \rightarrow
\left\langle \Lambda\left(  f\right)  ,\varphi\right\rangle
\]
is a bilinear map, separately continuous since $\Lambda$ is continuous. As
$\mathcal{S}\left(  \mathbb{R}^{n}\right)  $ and $\mathcal{S}\left(
\mathbb{R}^{m}\right)  $ are Fr\'{e}chet spaces, $\Theta$ is globally
continuous. There exist $C_{1}>0$,\ $\left(  q_{1},l_{1},q_{2},l_{2}\right)
\in\mathbb{N}^{4}$, such that
\[
\forall\left(  f,\varphi\right)  \in\mathcal{S}\left(  \mathbb{R}^{n}\right)
\times\mathcal{S}\left(  \mathbb{R}^{m}\right)  ,\;\;\;\left\vert \left\langle
\Lambda\left(  f\right)  ,\varphi\right\rangle \right\vert \leq C_{1}%
\,\mu_{q_{1},l_{1}}(f)\mu_{q_{2},l_{2}}(\varphi).
\]
In particular, for any $l\in\mathbb{N}$ and $\alpha\in\mathbb{N}^{m}$ with
$\left\vert \alpha\right\vert \leq l$, we have
\[
\forall x\in\mathbb{R}^{m},\;\left\vert \left\langle \Lambda\left(  f\right)
,\partial^{\alpha}\varphi_{\varepsilon}\left(  x-\cdot\right)  \right\rangle
\right\vert \leq C_{1}\,\mu_{q_{1},l_{1}}(f)\mu_{q_{2},l_{2}}(\partial
^{\alpha}\varphi_{\varepsilon^{s}}\left(  x-\cdot\right)  ),
\]
with
\begin{align*}
\left.  \forall x\in\mathbb{R}^{m},\;\mu_{q_{2},l_{2}}(\partial^{\alpha
}\varphi_{\varepsilon^{s}}\left(  x-\cdot\right)  )\right.   &  =\sup_{\xi
\in\mathbb{R}^{m},\left\vert \beta\right\vert \leq l_{2}}\left(  1+\left\vert
\xi\right\vert \right)  ^{q_{2}}\left\vert \partial^{\alpha+\beta}%
\varphi_{\varepsilon^{s}}\left(  x-\xi\right)  \right\vert \\
&  =\sup_{\xi\in\mathbb{R}^{m},\left\vert \beta\right\vert \leq l_{2}}\left(
1+\left\vert x-\xi\right\vert \right)  ^{q_{2}}\left\vert \partial
^{\alpha+\beta}\varphi_{\varepsilon^{s}}\left(  \xi\right)  \right\vert \\
&  \leq\left(  1+\left\vert x\right\vert \right)  ^{q_{2}}\mu_{q_{2},l_{2}%
+l}\left\vert \varphi_{\varepsilon^{s}}\right\vert .
\end{align*}
Using the definition of $\left(  \varphi_{\varepsilon^{s}}\right)
_{\varepsilon}$, we get $C_{2}>0$ such that $\mu_{q_{2},l_{2}+l}\left\vert
\varphi_{\varepsilon^{s}}\right\vert \leq C_{2}\,\varepsilon^{-s(m+l_{2}+l)}$,
for $\varepsilon$ small enough. Thus, there exists $C>0$, such that, for
$\varepsilon$ small enough,
\[
\forall x\in\mathbb{R}^{m},\ \ \ \left\vert \partial^{\alpha}\left(
L_{\varepsilon}(f)\right)  (x)\right\vert =\left\vert \left\langle
\Lambda\left(  f\right)  ,\partial^{\alpha}\varphi_{\varepsilon}\left(
x-\cdot\right)  \right\rangle \right\vert \leq C\left(  1+\left\vert
x\right\vert \right)  ^{q_{2}}\mu_{q_{1},l_{1}}(f)\varepsilon^{-s(m+l_{2}%
+l)}.
\]
Finally%
\[
\mu_{-q_{2},l}\left(  L_{\varepsilon}(f)\right)  \leq C\,\varepsilon
^{-s(m+l_{2}+l)}\mu_{q_{1},l_{1}}(f).
\]
The sequence $r\left(  \cdot\right)  =\left\{  l\mapsto s\left(
m+l_{2}+l\right)  \right\}  $ satisfies $\lim_{l\rightarrow+\infty}\left(
r(l)/l\right)  =s<1$. Recalling that $l_{1}$ does not depend on $l$, we obtain
our claim.$\,\medskip\smallskip$

\noindent\textit{Assertion (ii)}.-~We have the following:

\begin{lemma}
\label{GSTLmmGDequa}For all $u\in\mathcal{S}^{\prime}\left(  \mathbb{R}%
^{m}\right)  $, $\left[  \left(  u\ast\varphi_{\varepsilon^{s}}\right)
_{\varepsilon}\right]  $ is equal to $u$ in the generalized distribution sense.
\end{lemma}

\begin{proof}
Take $u\in\mathcal{S}^{\prime}\left(  \mathbb{R}^{m}\right)  $.\ There exist
$\alpha\in\mathbb{N}^{m}$, $q\in\mathbb{N}$, and $f:\mathbb{R}^{m}%
\rightarrow\mathbb{C}$ a continuous bounded function such that
\cite{Schwartz1}%
\[
u=\partial^{\alpha}\left(  M^{q}f\right)  ,
\]
where $M:\mathbb{R}^{m}\rightarrow\mathbb{C}$ is the function defined by
$M(x)=1+\left\vert x\right\vert ^{2}$. For any $g\in\mathcal{S}\left(
\mathbb{R}^{m}\right)  $, we have%
\[
\left\langle u\ast\varphi_{\varepsilon^{s}},g\right\rangle =\left\langle
u,g\ast\check{\varphi}_{\varepsilon^{s}}\right\rangle =\left(  -1\right)
^{\left\vert \alpha\right\vert }\left\langle f,M^{q}\,\left(  \partial
^{\alpha}g\right)  \ast\check{\varphi}_{\varepsilon^{s}}\right\rangle .
\]
On the other hand,
\[
\left\langle u,g\right\rangle =\left(  -1\right)  ^{\left\vert \alpha
\right\vert }\left\langle f,M^{q}\,\partial^{\alpha}g\right\rangle .
\]
Thus%
\[
\left\langle u\ast\varphi_{\varepsilon^{s}},g\right\rangle -\left\langle
u,g\right\rangle =\left\langle f,M^{q}\,\left(  \left(  \partial^{\alpha
}g\right)  \ast\check{\varphi}_{\varepsilon^{s}}-\partial^{\alpha}g\right)
\right\rangle .
\]
A simplification of the proof of \ref{LmnMTh4} shows that $\left(  \left(
\partial^{\alpha}g\right)  \ast\check{\varphi}_{\varepsilon^{s}}%
-\partial^{\alpha}g\right)  _{\varepsilon}\in\mathcal{N}_{\mathcal{S}}\left(
\mathbb{R}^{d}\right)  $.\ The same holds for $M^{q}\,\left(  \left(
\partial^{\alpha}g\right)  \ast\check{\varphi}_{\varepsilon^{s}}%
-\partial^{\alpha}g\right)  $. Thus, for all $p$ in $\mathbb{N}$,
\[
\left\langle u\ast\varphi_{\varepsilon^{s}},g\right\rangle -\left\langle
u,g\right\rangle =\mathrm{O}\left(  \varepsilon^{p}\right)  \text{ as
}\varepsilon\rightarrow0\text{.}%
\]
\medskip
\end{proof}

This lemma implies that for all $f\in\mathcal{S}\left(  \mathbb{R}^{n}\right)
$, $\left[  \left(  L_{\varepsilon}(f)\right)  _{\varepsilon}\right]  _{\tau
}=\left[  \left(  \Lambda\left(  f\right)  \ast\varphi_{\varepsilon^{s}%
}\right)  _{\varepsilon}\right]  _{\tau}$ is equal to $\Lambda\left(
f\right)  $ in the generalized distribution sense. On the other hand,
according to Theorem \ref{KernSTHKern1}, $\left[  \left(  L_{\varepsilon
}(f)\right)  _{\varepsilon}\right]  _{\tau}=\widetilde{H}_{L}\left(  f\right)
$ where $\widetilde{H}_{L}$ is the integral operator associated to the
canonical extension of $\left(  L_{\varepsilon}\right)  _{\varepsilon}$.\ This
ends the proof of Proposition \ref{GSTCoroTHS}.


\begin{thebibliography}{99}                                                                                               %


\bibitem {AntRad}\textsc{Antonevich A.B., Radyno Ya.V.} \emph{On a general
method for constructing algebras of generalized functions}\textit{.} Soviet
Math. Dokl. \textbf{43}(3):680-684, 1992.

\bibitem {BCD}\textsc{S.~Bernard, J.-F.~Colombeau, A.~Delcroix.}
\emph{Generalized Integral Operators and Applications}. Accepted in Math.
Proc. Cambridge Philos. Soc.

\bibitem {Col1}\textsc{J.-F.~Colombeau.} \emph{New Generalized Functions and
Multiplication of Distributions}. North-Holland (Amsterdam, Oxford, New-York), 1984.

\bibitem {Col2}\textsc{J.-F.~Colombeau.} \emph{Elementary introduction to New
generalized Functions.} North-Holland (Amsterdam, Oxford, New-York), 1985.

\bibitem {ADSKTT}\textsc{A.~Delcroix.} \emph{Generalized Integral Operators
and Schwartz Kernel Theorem}.{ }J. Math. Anal. Appl. \textbf{306}(2):481-501, 2005.

\bibitem {ADRegAppl}\textsc{A.~Delcroix.} \emph{Regular nonlinear generalized
functions and applications}. Accepted in Bull. Cl. Sci. Math. Nat. Sci.\ Math, 2005.

\bibitem {ADRapDec}\textsc{A.~Delcroix.} \emph{Regular rapidly decreasing
nonlinear generalized functions. Application to microlocal regularity}%
.\ Preprint AOC, 2005.

\bibitem {DHPV1}\textsc{A.~Delcroix, M.~Hasler, S.~Pilipovi\'{c},
V.~Valmorin.} \emph{Generalized function algebras as sequence space algebras}.
Proc. Amer. Math. Soc. \textbf{132}:2031-2038, 2004.

\bibitem {Garetto}\textsc{C.~Garetto.} \emph{Pseudo-differential Operators in
Algebras of Generalized Functions and Global Hypoellipticity}. Acta Appl.
Math. \textbf{80}(2):123-174, 2004.

\bibitem {GaGrOb}\textsc{C.~Garetto,~T.~Gramchev,~M.~Oberguggenberger.}
\emph{Pseudo-Differential operators and regularity theory}. Preprint 8-2003,
Preprint series of the department of Engineering Mathematics, Geometry and
Computer Science, University of Innsbruck, http://techmath.uibk.ac.at/mathematik/publikationen.

\bibitem {GKOS}\textsc{M.~Grosser, M.~Kunzinger, M.~Oberguggenberger,
R.~Steinbauer.} \emph{Geometric Theory of Generalized Functions with
Applications to General Relativity}. Kluwer Academic Press, 2001.

\bibitem {HorPDOT1}\textsc{L.~H\"{o}rmander.} \emph{The analysis of Linear
Partial Differential Operators I,\ distribution theory and Fourier Analysis}.
Grundlehren der mathematischen Wissenchaften 256. Springer Verlag (Berlin,
Heidelberg, New York), 2nd edition, 1990.

\bibitem {JAM2}\textsc{J.-A.~Marti.} \emph{Non linear Algebraic analysis of
delta shock wave to Burgers' equation}.\ Pacific\ J.\ Math.\ \textbf{210}%
(1):165-187, 2003.

\bibitem {NePiSc}\textsc{M.~Nedeljkov, S.~Pilipovi\'{c}, D.~Scarpal\'{e}zos.}%
\ \emph{The linear theory of Colombeau generalized functions.} Pitman Research
Notes in Mathematics Series 385.\ Longman, 1998\ 

\bibitem {Ober1}\textsc{M.~Oberguggenberger.} \emph{Multiplication of
Distributions and Applications to Partial Differential Equations.}\ Longman
Scientific \& Technical, 1992.

\bibitem {PScamb}\textsc{S.~Pilipovi\'{c}, D.~Scarpal\'{e}zos.}\ \emph{Colombeau
generalized ultradistributions}. Math.\ Proc.\ Cambridge\ Philos.\ Soc.
\textbf{130}:541-553, 2001.

\bibitem {Radyno}\textsc{Ya.-V. Radyno, Ngo Fu Tkhan, S.\ Ramadan}. \emph{The
Fourier transform in an algebra of new generalized functions}. Russian acad.
Sci. Dokl. Math. \textbf{46}(3):414-417, 1992

\bibitem {Scarpa1}\textsc{D.~Scarpal\'{e}zos.}\ \emph{Colombeau's generalized
functions: Topological structures; Microlocal properties.\ A simplified point
of view.} Pr\'{e}publication Math\'{e}matiques de Paris 7/CNRS, URA212, 1993.

\bibitem {Scarpa2}\textsc{D.~Scarpal\'{e}zos.}\ \emph{Colombeau's generalized
functions: Topological structures; Microlocal properties.\ A simplified point
of view. Part I}. Bull. Cl. Sci. Math. Nat. Sci.\ Math. 25:89-114, 2000.

\bibitem {Scarpa3}\textsc{D.~Scarpal\'{e}zos.}\emph{\ Colombeau's generalized
functions: Topological structures; Microlocal properties.\ A simplified point
of view. Part II}. Publ. Inst.\ Math. (Beograd) (N.S.) \textbf{76}%
(90):111-125, 2004.

\bibitem {Schwartz1}\textsc{L.~Schwartz.} \emph{Th\'{e}orie des Distributions}%
.\ Hermann (Paris), 3rd print, 1965.{}

\bibitem {Treves1}\textsc{F.~Treves.} \emph{Topological vector spaces,
distributions and kernels}. Academic Press (New York, London), 1967.

\bibitem {Valmo1}\textsc{V.~Valmorin.} \emph{Schwartz kernel theorem in
Colombeau algebras}.\ Preprint AOC, 2005.
\end{thebibliography}
\end{document}